\documentstyle[12pt]{article}
\begin{document}
\title{Graphs and the (co)homology of Lie algebras}
\author {Qibing Zheng\\Department of Mathematics, Nankai University\\
Tianjin 300071, China\\
zhengqb@nankai.edu.cn\footnote{Project Supported by Natural Science
Foundation of China}} \maketitle
\input amssym.def
\newsymbol\leqslant 1336
\newsymbol\geqslant 133E
\baselineskip=20pt

\begin{abstract} In this paper, we develop a diamond graph theory
and apply the theory to the (co)homology of the Lie algebra
generated by positive systems of the classical semi-simple Lie
algebras over the field of complex numbers. As an application, we
give the weight decomposition of  the diamond Lie algebra with
Dynkin graph $A_{n+1}$ and compute the rank of every weight subgraph
of it.
\end{abstract}

\newtheorem{Definition}{Definition}[section]
\newtheorem{Theorem}{Theorem}[section]
\newtheorem{Lemma}{Lemma}[section]
\newtheorem{Example}{Example}[section]
\font\hua=eusm10 scaled\magstephalf

The classification of compact simply connected Lie groups over $\Bbb
C$ is due to the classification of the Dynkin graph of the
associated semi-simple Lie algebra. The Dynkin graph determines a
unique positive root system and Kostant in \cite{k} tells that the
integral cohomology of the Lie subalgebra generated by the positive
root system is in 1-1 correspondence with the group ring of the Weyl
group of the Lie algebra. But the torsion part of the (co)homology
of the Lie algebra generated by a positive root system is also very
important. For example, let $\frak A_n$ be the Lie algebra generated
by the positive root system with Dynkin graph $A_n$ and $\frak
A_{\infty}=\cup_n \frak A_n$ (a graded Lie algebra). Then,
$H_*(\frak A_{n};\Bbb Z_p)$ is the $E_2$-term of the spectral
sequence induced by the lower central series converging to the
homology of the group of integral upper-triangular matrices in
\cite{v} and $H^*(\frak A_{\infty};\Bbb Z_p)$ ($\Bbb Z_p$ the field
of integers modular a prime $p$) is a direct sum summand of the
$E_1$-term of May spectral sequence in \cite{m} converging to the
cohomology of the Steenrod algebra.

In this paper, we develop a diamond graph theory and study the
(co)homology of the Lie algebra generated by a diamond root system.
Diamond graphs give more information about the torsion part the the
(co)homology of these Lie algebras. The main results of this paper
is Theorem 2.2 and Theorem 3.5.\vspace{4mm}

\section{Graphs and (co)chain complexes}\vspace{4mm}

\hspace{6mm}In this paper, all objects are of finite type. Graphs
are finite graphs and Abellian groups are finitely generated.
$\otimes$ means $\otimes_{\Bbb Z}$ with $\Bbb Z$ the ring of
integers. For a set $S$, $\Bbb Z(S)$ denotes the free Abellian group
generated by $S$. $S$ is called a base of $\Bbb Z(S)$ .\vspace{4mm}

Recall that a graph is a pair of sets $(G,E)$ such that every
element of $E$ is a subset of two elements of $G$. The elements of
$G$ are called vertices of the graph and the elements of $E$ are
called edges of the graph. We omit the set of edges and simply
denote by $G$ the graph $(G,E)$. We denote by $[a,b]$ an edge and
call $a$ and $b$ neighbors of each other. A graph is finite if the
set of vertices is a finite set. \vspace{4mm}

\begin{Definition}{\rm For a graph $G$, a gradation $|\cdot|$ on
$G$ is a map $|\cdot|\colon G\to \Bbb Z$ such that if $[a,b]$ is an
edge, then $|a|{-}|b|=\pm 1$. A graph is called gradable if there is
a gradation on it. A graded graph with base graph $G$ is a pair
$(G,|\cdot|)$ with $G$ a gradable graph and $|\cdot|$ a gradation on
it. $(G,|\cdot|)$ is simply denoted by $|G|$. Two graded graphs are
isomorphic if there is a graph isomorphism $f$ that keeps the
gradation, i.e., $|f(a)|-|f(b)|=|a|-|b|$ for all $a,b$. A vertex $v$
of a graded graph $|G|$ is called a top vertex if $|u|=|v|-1$ for
every neighbor $u$ of $v$. A vertex $v$ of $|G|$ is called a bottom
vertex if $|u|=|v|+1$ for every neighbor $u$ of
$v$.}\end{Definition}\vspace{4mm}

Notice that a gradable graph has no triangle as a subgraph. In fact,
it has no $(2n{+}1)$-polygon as a subgraph. That is, there is no
$2n{+}1$ vertices $v_0,v_1,\cdots,v_{2n}$ such that
$[v_0,v_{1}],\cdots,$\\$[v_{2n-1},v_{2n}],[v_{2n},v_{0}]$ are all
edges.\vspace{4mm}

\begin{Definition}{\rm
Let $|G|$ be a graded graph. If $v$ is a top vertex of $G$, then
$|\cdot|_1$ defined by $|v|_1=|v|{-}2$ and $|u|_1=|u|$ for every
other vertex $u$ is also a gradation on $G$ which is called the
lowering of $|\cdot|$ by the vertex $v$. The graded graph $|G|_1$ is
called the lowering of $|G|$ by the vertex $v$. If $v$ is a bottom
vertex of $|G|$, then $|\cdot|_2$ defined by $|v|_2=|v|{+}2$ and
$|u|_2=|u|$ for every other vertex $u$ is also a gradation on $G$
which is called the lifting of $|\cdot|$ by the vertex $v$. The
graded graph $|G|_2$ is called the lifting of $|G|$ by the vertex
$v$. Two gradations on $G$ are equivalent if one of them can be
obtained from the other by a finite composite of lowerings and
liftings. Two graded graphs are equivalent if there is a graph
isomorphism that induces a gradation
equivalence.}\end{Definition}\vspace{4mm}

\begin{Theorem}{\rm Let $G$ be a gradable graph.
Two gradations $|\cdot|$ and $|\cdot|'$ on $G$ are equivalent if and
only if $|v|-|v|'$ is even for all $v\in G$. }\end{Theorem}

Proof. The necessary part is by definition. Now we prove that if
$|w|-|w|'$ is even for all vertices $w$, then $|\cdot|$ and
$|\cdot|'$ are equivalent. Let $b_1$ be a vertex such that
$|b_1|\geqslant |w|$ for all $w\in G$. If $|b_1|>1$, then lower
$|\cdot|$ by $b_1$ and we get a new gradation $|\cdot|_1$. Let $b_2$
be a vertex such that $|b_2|_1\geqslant |w|_1$ for all $w\in G$. If
$|b_2|_1>1$, then lower $|\cdot|_1$ by $b_2$ and we get a new
gradation $|\cdot|_2$. Repeat this process if there is vertex with
degree $>1$. Since $G$ is a finite graph, this process will come to
an end. That is, there exists an $n$ and gradations $|\cdot|_1$,
$\cdots$, $|\cdot|_n$ such that each $|\cdot|_{i+1}$ is a lowering
of $|\cdot|_i$ by the vertex $b_{i+1}$ and for all $v\in G$,
$|v|_n\leqslant 1$. Similarly, by lifting the smallest degree
vertex, we get an $m$ and gradations $|\cdot|_{n+1}$, $\cdots$,
$|\cdot|_{n+m}$ such that each $|G|_{n+i+1}$ is a lifting of
$|G|_{n+i}$ by the vertex $b_{n+i+1}$ and for all $w\in G$,
$|w|_{n+m}=0$ or $1$. $|\cdot|$ is equivalent to $|\cdot|_{n+m}$.
Similarly, $|\cdot|'$ is equivalent to a gradation $|\cdot|'_{s+t}$
such that $|w|'_{s+t}=0$ or $1$ for all $w\in G$. Since $|w|-|w|'$
is even for all $w\in G$, we have $|\cdot|_{n+m}=|\cdot|'_{s+t}$.
Thus, $|\cdot|$ is equivalent to $|\cdot|'$. \hfill
Q.E.D.\vspace{4mm}

Recall that a path from $a$ to $b$ is a sequence of vertices
$a=v_0,v_1,\cdots,v_{n-1},v_n=b$ such that either $[v_{i-1},v_{i}]$
is an edge, or $v_{i-1}=v_{i}$ for $i=1,\cdots,n$. The length of the
path is the number of edges $[v_{i-1},v_i]$. The distance $d(a,b)$
between two vertices $a$ and $b$ is the minimum of lengths of all
paths from $a$ to $b$. If there is no path from $a$ to $b$, we
define $d(a,b)=\infty$. A graph is connected if the distance between
every pair of its vertices is finite.

\begin{Theorem}{\rm Let $G$ be a connected graph with more
than one vertex. $G$ is gradable if and only if its vertex set has a
unique distance decomposition $G=G_1\sqcup G_2$ ($\sqcup$ is the
disjoint union) such that for all $u,v\in G_i$, $i=1,2$, $d(u,v)$ is
even and for all $a\in G_1$ and $b\in G_2$, $d(a,b)$ is odd. $G_1$
and $G_2$ are called the distance components of $G$.}\end{Theorem}

Proof. If there is a distance decomposition $G=G_1\sqcup G_2$, then
the gradation $|\cdot|$ defined by $|a|=0$ for all $a\in G_1$ and
$|b|=1$ for all $b\in G_2$ is a gradation. So $G$ is gradable.

If $G$ is gradable, then from the proof of Theorem 1.1 we have that
there are only two equivalent classes of gradations on $G$
represented by the two gradations $|\cdot|_1$ and $|\cdot|_2$ such
that $|v|_1{+}|v|_2=1$ and $|v|_i=0$ or $1$ for all $v\in G$. Then
$G_1=\{v{\in}G\,|\,|v|_1=0\}$ and $G_2=\{v{\in}G\,|\,|v|_1=1\}$ are
the distance components. \hfill Q.E.D.\vspace{4mm}

\begin{Definition}{\rm
For a gradable graph $G$, any gradation $|\cdot|$ satisfying that
$|v|=0$ or $1$ for all $v\in G$ is called a representation gradation
of $G$. The set $G_1=\{v\in G\,|\,|v|=0\}$ and $G_2=\{v\in
G\,|\,|v|=1\}$ are called the distance components of the
representation gradation.
 }\end{Definition}\vspace{4mm}

\begin{Theorem}{\rm Let $|G|$ be a connected graded graph.
If $|G|$ has only one bottom vertex $v$, then $d(u,v)=|u|-|v|$ for
all vertex $u$. Such a graded graph is called a positive distance
graph relative to $v$. If $|G|$ has only one top vertex $v$, then
$d(u,v)=|v|-|u|$ for all vertex $u$. Such a graded graph is called a
negative distance graph relative to $v$.}\end{Theorem}

Proof. We only prove the positive distance case. Suppose $|\cdot|$
is a gradation that has only one bottom vertex $v$. Let $N={\rm
min}\{|u|\,|\,u\in G\}$. If $|u|=N$, then $u$ is a bottom vertex and
so $u=v$. This implies that for all $u\in G$, $|u|\geqslant |v|$ and
the equality holds if and only if $u=v$. We use induction on $n$ to
prove that $|u|=|v|{+}n$ if and only if $d(u,v)=n$. If $n=0,1$, the
conclusion is trivial. Suppose for some $n>1$, we have
$|u'|=|v|{+}i$ if and only if $d(u',v)=i$ for $i=0,1,\cdots,n$. Then
for $|u|=|v|{+}n{+}1$, the induction hypothesis implies that
$d(u,v)>n$. Since $u$ is not a bottom vertex, there is a neighbor
$w$ of $u$ such that $|w|=|v|{+}n$. By the induction hypothesis,
$d(w,v)=n$. So $d(u,v)\leqslant d(u,w){+}d(w,v)=n{+}1$. Thus,
$d(u,v)=n{+}1$. The conclusion holds. \hfill Q.E.D.\vspace{4mm}

\begin{Definition}{\rm
Let $G$ be a gradable graph. A connection $\nu$ on $G$ is a map
$\nu\,\colon G\times G\to\Bbb Z$ that satisfies the following
conditions.

1) $\nu(a,b)=\nu(b,a)$ for all $a,b\in G$.

2) $\nu(a,b)\neq 0$ if and only if $[a,b]$ is an edge of $G$.

Two connections $\nu,\nu'$ are equivalent if there is a map $e\colon
G\to\{\pm 1\}$ such that $\nu(a,b)=e(a)e(b)\nu'(a,b)$ for all
$a,b\in G$.

A graph with connection is a pair $(G,\nu)$ with $G$ a gradable
graph and $\nu$ a connection on $G$. }\end{Definition}\vspace{4mm}

Notice that a connection can be defined on an ungradable graph. But
such a map has no representation matrix defined as follows. So we
define connection only on gradable graphs. \vspace{4mm}

\begin{Definition}{\rm
For a graph with connection $(G,\nu)$, its representation matrix
$A=(a_{i,j})_{m\times n}$ is defined as follows. If $G$ has only one
vertex, $A=(0)_{1\times 1}$, the $1{\times} 1$ zero matrix. The
global dimension $D(A)$ of $A$ is defined to be $1$. If $G$ is
connected and has more than one vertices, take a representation
gradation $|\cdot|$ of $G$ and suppose $v_1,\cdots,v_m$ and
$w_1,\cdots,w_n$ are the distance components of $|\cdot|$, then
$a_{i,j}=\nu(v_i,w_j)$. The global dimension $D(A)$ of $A$ is
defined to be $m{+}n{-}2r$, where $r$ is the rank of $A$. If $G$ is
not connected, then its representation matrix is the direct sum (see
the next definition) of all its connected component representation
matrices and the global dimension of $A$ is the sum of the global
dimensions of all its connected component representation
matrices.}\end{Definition}\vspace{4mm}

The representation matrix is not unique and depends on the order of
distance components and their vertices. Different equivalent
connections have different representation matrices. To make the
representation matrices unique under equivalences, we have the
following definition.

\begin{Definition}{\rm
Two matrices over $\Bbb Z$ are equivalent if one of them can be
obtained from the other by a finite composite of the following
transformations.

1) Permute the rows of the matrix.

2) Replace a row $\alpha$ of the matrix by $-\alpha$.

3) Replace the matrix $A$ by its transpose matrix $A^T$.

4) Replace matrix $\left(\begin{array}{cc}A&0\\0&B\end{array}\right)$ by $\left(\begin{array}{cc}A^T&0\\0&B\end{array}\right)$.

For two  matrices $A=(a_{i,j})_{m_1{\times}n_1}$ and
$B=(b_{k,l})_{m_2{\times}n_2}$, their direct sum is the matrix
\vspace{-2mm}
 \hspace*{60mm}$$A\oplus
B=\left(\begin{array}{cc}A&0\\0&B\end{array}\right)_{(m_1+m_2)\times(n_1+n_2)}$$
and their orthogonal product $A{\times}B=(c_{s,t})_{m{\times}n}$ is
the following matrix. $m=m_1m_2{+}n_1n_2$, $n=m_1n_2{+}n_1m_2$,
$c_{s,t}=0$ except the following,\vspace{-2mm}
\begin{eqnarray*}&&c_{im_1+k,jn_1+k}=a_{i,j},\hspace{22.5mm}\,i=1,\cdots,m_1,\,j=1,\cdots,n_1,\,k=1,\cdots,n_2\\
&&c_{m_1n_2+jn_1+k,m_1m_2+im_1+k}=a_{i,j},\hspace{2mm}\,i=1,\cdots,m_1,\,j=1,\cdots,n_1,\,k=1,\cdots,m_2\\
&&c_{m_1n_2+kn_1+i,kn_1+j}=b_{i,j},\hspace{14mm}\,i=1,\cdots,m_2,\,j=1,\cdots,n_2,\,k=1,\cdots,n_1\\
&&c_{km_1+j,m_1m_2+km_1+i}=-b_{i,j},\hspace{8mm}\,i=1,\cdots,m_2,\,j=1,\cdots,n_2,\,k=1,\cdots,m_1
\end{eqnarray*}
 }\end{Definition}\vspace{4mm}

It is obvious that if $A$ and $B$ are respectively equivalent to
$A'$ and $B'$, $A\oplus B$ is equivalent to $A'\oplus B'$ and
$A\times B$ is equivalent to $A'\times B'$. \vspace{4mm}

\begin{Definition}{\rm A chain graph is a pair $(|G|,\nu)$ with
$|G|$ a graded graph and $\nu$ a connection on $G$ such that for
every pair $a,b\in G$ with $|b|=|a|{+}2$,
$\sum_{|c|=|a|+1}\nu(a,c)\nu(c,b)=0$. Two chain graphs
$(|G_1|_1,\nu_1)$ and $(|G_2|_2,\nu_2)$ are similar if there is a
graph isomorphism $\psi\colon G_1\to G_2$ such that the induced map
$\nu'_1(a,b)=\nu_2(\psi(a),\psi(b))$ for all $a,b\in G_1$ is a
connection equivalent to $\nu_1$. If $\psi$ is a graded graph
isomorphism, $(|G_1|_1,\nu_1)$ and $(|G_2|_2,\nu_2)$ are
isomorphic.}\end{Definition}\vspace{4mm}

\begin{Definition}{\rm
For a chain graph $(|G|,\nu)$, the associated chain complex $(\Bbb
Z(G),d)$  and associated cochain complex $(\Bbb Z(G),\delta)$ of the
graph is defined as follows. $dv=\sum_{|w|=|v|-1}\nu(v,w)w$ and
$\delta v=\sum_{|w|=|v|+1}\nu(v,w)w$ for all $v\in G$.

The homology of $(\Bbb Z(G)\otimes R,d)$ with $R$ a commutative ring
is called the homology of $(|G|,\nu)$ over the coefficient ring $R$
and is denoted by $H_*(|G|;R)$ and $H_*(|G|)=H_*(|G|;\Bbb Z)$.
Dually, $H^*(|G|;R)=H^*({\rm Hom}_{\Bbb Z}(\Bbb Z(G),R),\delta)$ is
the cohomology of $(|G|,\nu)$ over the coefficient ring $R$ and
$H^*(|G|)=H^*(\Bbb Z(G),\delta)=H^*(|G|;\Bbb Z)$.
}\end{Definition}\vspace{4mm}

\begin{Definition}{\rm
Let $(|G_1|_1,\nu_1)$ and $(|G_2|_2,\nu_2)$ be two chain graphs. The
disjoint union graph $(|G_1\sqcup G_2|,\nu)$ is defined as follows.
The restriction of $|\cdot|$ on $G_i$ is $|\cdot|_i$ and the
restriction of $\nu$ on $G_i\times G_i$ is $\nu_i$ and $\nu(a,b)=0$
for all $a\in G_1$ and $b\in G_2$. The product graph $(|G_1\times
G_2|,\nu)$ is defined as follows. $|(g_1,g_2)|=|g_1|_1{+}|g_2|_2$
for all $g_i\in G_i$. $[(g_1,g_2),(g'_1,g'_2)]$ is an edge of
$G_1\times G_2$ if either $[g_1,g'_1]$ is an edge $G$ and
$g_2=g'_2$, or $g_1=g'_1$, $[g_2,g'_2]$ is an edge $G_2$. For all
$g_i,g'_i\in G_i$, $\nu((g_1,g_2),(g'_1,g_2))=\nu_1(g_1,g'_1)$,
$\nu((g_1,g_2),(g_1,g'_2))=(-1)^{|g_1|_1}\nu_2(g_2,g'_2)$.
}\end{Definition}\vspace{4mm}

\begin{Theorem}{\rm
Let $(|G|,\nu)$ and $(|H|,\nu)$ be two chain graphs. The associated
chain complex $(\Bbb Z(G{\sqcup}H),d)$ of the disjoint union graph
$(|G{\sqcup}H|,\nu)$ is the direct sum chain complex $(\Bbb
Z(G)\oplus\Bbb Z(H),d)$. The associated chain complex $(\Bbb
Z(G{\times}H),d)$ of the product graph $(|G{\times}H|,\nu)$ is the
tensor product chain complex $(\Bbb Z(G)\otimes\Bbb Z(H),d)$. The
same conclusion holds for associated cochain
complexes.}\end{Theorem}

Proof. Direct checking.\hfill Q.E.D.\vspace{4mm}

\begin{Theorem}{\rm Let $(|G|,\nu)$ be a chain graph. Then \\
\hspace*{24mm}$D(|G|,\nu)=\sum_{k=-\infty}^{\infty}{\rm
dim}H_k(|G|)=\sum_{k=-\infty}^{\infty}{\rm dim}H^k(|G|)=D(A)$\\for
any representation matrix $A$ of the graph with connection
$(G,\nu)$, where dim$H$ means the dimension of the free part of the
(co)homology. Thus, the global dimension $D(G,\nu)=D(|G|,\nu)$ of
the graph with connection $(G,\nu)$ is well-defined and for two
similar chain graphs $(|G_1|_1,\nu_1)$ and $(|G_2|_2,\nu_2)$,
$D(|G_1|_1,\nu_1)=D(|G_2|_2,\nu_2)=D(G_1,\nu_1)$. }\end{Theorem}

Proof. We may suppose $G$ is connected and only prove the homology
case. If $|\cdot|$ is a representation gradation, the conclusion is
obvious. Suppose $|\cdot|$ is not a representation gradation and we
may suppose the associated chain complex $(C,d)$ of the chain graph
$(G,\nu)$ is the following.

\[0\to C_n\stackrel{d_n}{\longrightarrow} C_{n-1}\stackrel{d_{n-1}\,\,\,\,}{\longrightarrow}
 C_{n-2}\stackrel{d_{n-2}\,\,\,\,}{\longrightarrow}\cdots\stackrel{d_0}{\to} C_0\to 0\]

Define chain complex $(C',d)$

\[0\to C_{n-1}\stackrel{d'_{n-1}\,\,\,\,}{\longrightarrow}
 C_{n-2}{\oplus}C_n\stackrel{d'_{n-2}\,\,\,\,}{\longrightarrow}\cdots\stackrel{d'_0}{\to} C_0\to 0\]
as follows. Suppose $C_n$ has a base $a_1,\cdots,a_s$ and $C_{n-1}$
has a base $b_1,\cdots,b_t$ and $d_n(a_i)=\sum_j c_{i,j}b_j$. Then
$d'_{n-1}(a_i)=0$ and $d'_{n-1}(b_i)=d_{n-1}(b_i){+}\sum_k
c_{k,i}a_k$ and $d'_{k}=d_k$ for $k<n{-}2$. It is obvious that
$(C',d)$ is the associated chain complex of another chain graph
$(|G|',\nu)$. Since $D(|G|)=\sum_{k=-\infty}^{\infty}{\rm
dim}H_k(|G|;F)$ for all field $F$ of characteristic $0$, we compute
the number for $F=\Bbb R$, the real number. In this case, we may
suppose all the $C_n$ are vector spaces over $\Bbb R$. Let
$\delta_n=d'_{n-1}{-}d_{n-1}$ and $A=(c_{i,j})_{s\times t}$. $AA^T$
is a symmetric matrix and there is an orthogonal matrix $Q$ such
that $AA^T=Q^TDQ$ with $D=diag(k_1,\cdots,k_r,0,\cdots,0)$ a
diagonal matrix ($k_i>0$). Suppose $Q=(q_{i,j})_{s\times s}$ and let
$a'_i=\sum_j q_{i,j}a_j$. Then $a'_1,\cdots,a'_s$ is another base of
$C_n$ such that if $d_n(a'_i)=\sum_j c'_{i,j}b_j$, then
$\delta_n(b_i)=\sum_k c'_{k,i}a'_k$ and that
$\delta_nd_n(a'_i)=k_ia'_i$ for $i=1,\cdots,r$,
$\delta_nd_n(a'_j)=0$ for $j=r{+}1,\cdots,s$. By definition,
$d_n(a'_j)=0$ for $j=r{+}1,\cdots,s$. So dim$H_n(C;\Bbb R)=s{-}r$,
dim$H_n(C';\Bbb R)=0$, dim$H_{n-1}(C;\Bbb R)=$ dim$H_{n-1}(C';\Bbb
R)$, dim$H_{n-2}(C';\Bbb R)=$ dim$H_{n-2}(C;\Bbb R){+}s{-}r$,
dim$H_{k}(C;\Bbb R)=$ dim$H_{k}(C';\Bbb R)$ otherwise. Thus, $(C,d)$
and $(C',d)$ have the same global dimension. Repeat the above
process to $(C',d)$ and we can prove that $(C,d)$ and the associated
chain complex of the representation gradation have the same global
dimension.
 \hfill Q.E.D.\vspace{4mm}

\begin{Definition}{\rm
Let $G$ be a gradable graph. A connection $\nu$ on $G$ is deformable
if it satisfies that for all $a,b\in G$ with $d(a,b)=2$, $\sum_{c\in
G}\nu(a,c)\nu(c,b)=0$.

A gradable graph is deformable if there is a deformable connection
on it. A deformation graph is a pair $(G,\nu)$ with $G$ a gradable
graph and $\nu$ a deformable connection on $G$. Two deformation
graphs are isomorphic if there is a graph isomorphism that induces a
deformable connection equivalence. }\end{Definition}\vspace{4mm}

\begin{Theorem}{\rm For a deformation graph $(G,\nu)$, all its
representation matrices are equivalent $n{\times}n$ square matrices.
If $\nu'$ is a deformable connection equivalent to $\nu$, all the
representation matrices of $(G,\nu')$ are equivalent to that of
$(G,\nu)$. The equivalent matrix class is called the representation
class of $(G,\nu)$. The representation class of the disjoint union
of two graphs is the direct sum class of the two graphs. The
representation class of the product of two graphs is the orthogonal
product class of the two graphs. }\end{Theorem}

Proof. We only prove that representation matrices are square
matrices. Other conclusions are direct checkings. We may suppose $G$
is connected with more than one vertices. $\nu$ is a deformable
connection implies that every representation matrix $A$ satisfies
that $AA^T$ and $A^TA$ are both diagonal matrices with positive
diagonal entries. This implies that $A$ is a square matrix.

\hfill Q.E.D.\vspace{4mm}

\begin{Definition}{\rm The rank $r(v)$ of a vertex $v$ of
a deformation graph $(G,\nu)$ is the non-negative integer
$r(v)=\sum_{w\in G}\nu(w,v)^2$. }\end{Definition}\vspace{4mm}

\begin{Theorem}{\rm For a connected deformation graph $(G,\nu)$, all its
vertex $v$ have the same rank which is called the rank of $(G,\nu)$
and is denoted by $r(G,\nu)$.

 }\end{Theorem}

Proof. If $G$ has only one vertex, then by definition $r(G,\nu)=0$.
Suppose $G$ has more than one vertices and $v_1,\cdots,v_n$ and
$w_1,\cdots,w_n$ are the distance components of $G$ with
representation matrix $A=(\,a_{i,j}{=}\nu(v_i,w_j)\,)_{n\times n}$.
Then $AA^T={\rm diag}(d_1,\cdots,d_n)$ and $A^TA={\rm
diag}(d'_1,\cdots,d'_n)$ with $d_i=r(v_i)$ and $d'_j=r(w_j)$, where
diag$(\cdots)$ represents the diagonal square matrix. Thus,
$AA^TA={\rm diag}(d_1,\cdots,d_n)A=A{\rm diag}(d'_1,\cdots,d'_n)$
and $a_{i,j}\neq 0$ implies $d_i=d'_j$. Since $G$ is connected, for
any $i,j$, there is a path
$v_i,w_{j_1},v_{i_1},\cdots,w_{j_s},v_{i_s},w_j$ from $v_i$ to
$w_j$. So
$a_{i,j_1},a_{i_1,j_1},a_{i_1,j_2},a_{i_2,j_2},\cdots,a_{i_{s},j_s},a_{i_s,j}$
are all non-zero and
$d_i=d'_{j_1}=d_{i_1}=d'_{j_2}=\cdots=d'_{j_s}=d_{i_s}=d'_j$.
Similarly, there is a path from $v_i$ to $v_j$ for $i\neq j$ and we
have $d_i=d_j$. So $AA^T=A^TA=rE$ ($E$ unit matrix) and
$r=r(v_i)=r(w_j)$ for all $i,j$. \hfill Q.E.D.\vspace{4mm}

\begin{Theorem}{\rm Let $(G,\nu)$ be a deformation graph. Then for any gradation
$|\cdot|$ on $G$, $(|G|,\nu)$ is a chain graph. Such a chain graph
is called a deformable chain graph. }\end{Theorem}

Proof. $\nu$ is a deformable connection implies that for any
gradation $|\cdot|$ on $G$ and $a,b\in G$ with $|b|-|a|=\pm 2$,
$\sum_{|c|=|b|+1}\nu(a,c)\nu(c,b)=\sum_{c\in G}\nu(a,c)\nu(c,b)=0$.
\hfill Q.E.D.\vspace{4mm}

\begin{Theorem}{\rm For a connected deformable chain graph $(|G|,\nu)$
with rank $>0$, both $H_*(|G|)$ and $H^*(|G|)$ are torsion groups
and $H^k(|G|)=H_{k-1}(|G|)$ for all $k$.}\end{Theorem}

Proof. Since the representation matrix $A$ of $G$ is an orthogonal
matrix and so its global dimension $D(A)=0$. By Theorem 1.5,
$D(G,\nu)=0$ and so the free part of $H_k(|G|)$ and $H^k(|G|)$ are
trivial. By universal coefficient theorem, $H^k(|G|)={\rm
Ext}(H_{k-1}(|G|),\Bbb Z)=H_{k-1}(|G|)$ for all $k$. \hfill
Q.E.D.\vspace{4mm}

\begin{Theorem}{\rm Let $(|G|_1,\nu)$ be a connected deformable chain graph
with rank $n>0$. If $v$ is a bottom vertex of $|G|_1$ with $|v|_1=q$
and $|G|_2$ is the lifting of $|G|_1$ by $v$,  then there is a
divisor $k$ of $n$ such that ($\langle m\rangle$ denotes the group
of integers modular $m$)
\begin{eqnarray*}&&H_q(|G|_1)/\langle k\rangle=H_q(|G|_2),\\
&&H_{q+1}(|G|_2)/H_{q+1}(|G|_1)=\langle n/k\rangle,\\
&&H_i(|G|_1)=H_i(|G|_2)\,\,{\rm if}\,\,i\neq q,q{+}1.
\end{eqnarray*}

Dually, if $v$ is a top vertex with $|v|_1=q$ and $|G|_2$ is the
lowering of $|G|_1$ by $v$, then there is a divisor $k$ of $n$ such
that
\begin{eqnarray*}&&H^q(|G|_1)/\langle k\rangle=H^q(|G|_2),\\
&&H^{q-1}(|G|_2)/H^{q-1}(|G|_1)=\langle n/k\rangle,\\
&&H^i(|G|_1)=H^i(|G|_2)\,\,{\rm if}\,\,i\neq q,q{-}1.\end{eqnarray*}
 }\end{Theorem}

Proof. We only prove the cohomology case. Suppose
$(C_i,\delta_i,|\cdot|_i)$ is the associated cochain complex of
$(|G|_i,\nu)$. It is obvious that the free Abellian group generated
by $v$ is a cochain subcomplex of $C_1$ and that the free Abellian
group generated by vertices other than $v$ is a cochain subcomplex
of $C_2$. Denote the first cochain subcomplex of $C_1$ by $T_1$,
then the quotient cochain complex $C_1/T_1$ is just the second
cochain subcomplex of $C_2$ generated by vertices other than $v$.
Denote $C_1/T_1$ by $\tilde C$ and $C/\tilde C$ by $T_2$. Notice
that $H^{q-2}(T_2)=H^{q}(T_1)=\Bbb Z$; $H^{s}(T_2)=H^{t}(T_1)=0$,
otherwise. From the two short exact sequences of cochain complexes
$0\to\tilde C\to C_2\to T_2\to 0$ and $0\to T_1\to C_1\to\tilde C\to
0$ , we have two exact sequences of Abellian groups
\begin{eqnarray*}
0\to H^{q-2}(\tilde C)\to H^{q-2}(C_2)\to &H^{q-2}(T_2)&
\stackrel{\tau}{\to} H^{q-1}(\tilde C)\to H^{q-1}(C_2)\to 0\\
0\to H^{q{-}1}(C_1)\to H^{q{-}1}(\tilde C)
\stackrel{\pi}{\to}&H^{q}(T_1)&\to H^{q}(C_1)\to H^{q}(\tilde C)\to
0
\end{eqnarray*}
and have that $H^s(C_2)=H^s(\tilde C)$ if $s\neq q{-}1,q{-}2$ and
that $H^s(C_1)=H^s(\tilde C)$ if $s\neq q{-}1,q$. So the two exact
sequence are
\begin{eqnarray*}
0\to H^{q-2}(C_1)\to H^{q-2}(C_2)\to &\Bbb Z&
\stackrel{\tau}{\to} H^{q-1}(\tilde C)\to H^{q-1}(C_2)\to 0\quad (*)\\
0\to H^{q{-}1}(C_1)\to H^{q{-}1}(\tilde C) \stackrel{\pi}{\to}&\Bbb
Z&\to H^{q}(C_1)\to H^{q}(C_2)\to 0\quad (**)
\end{eqnarray*}
Since $H^q(C_1)$ is a torsion group, we have that $\pi$ is an
epimorphism. Thus, there is an integer $k$ such that $(**)$ implies
two short exact sequences
\begin{eqnarray*}
0\to H^{q{-}1}(C_1)\to H^{q{-}1}(\tilde C)\stackrel{\pi}{\to}\Bbb
Z\to 0, \quad 0\to\langle k\rangle\to H^{q}(C_1)\to H^{q}(C_2)\to 0
\end{eqnarray*}

Since $H^{q-1}(C_1)$ is a torsion group, we have that
$H^{q{-}1}(\tilde C)=H^{q{-}1}(C_1)\oplus\Bbb Z$ and the generator
$x$ of the free subgroup of $H^{q{-}1}(\tilde C)$ satisfies
$\pi(x)=k[v]$, where $[v]$ denote the cohomology class represented
by $v$.

Since $H^{q-2}(C_2)$ is a torsion group, $(*)$ implies two short
exact sequences
\begin{eqnarray*}
0\to H^{q-2}(C_1)\to H^{q-2}(C_2)\to 0,\quad 0\to\Bbb
Z\stackrel{\tau}{\to} H^{q-1}(C_1)\oplus\Bbb Z\to H^{q-1}(C_2)\to 0
\end{eqnarray*}

Since $H^{q-1}(C_2)$ is a torsion group, we have that $H^{q-1}(C_1)$
is a subgroup of $H^{q-1}(C_2)$ and
$H^{q-1}(C_2)/H^{q-1}(C_1)=(H^{q-1}(C_1)\oplus\Bbb Z)/(H^{q-1}(C_1)+
{\rm im}\tau)=\langle k'\rangle$ and $k'$ satisfies that
$\tau([v])\equiv k'x$ mod$H^{q-1}(C_1)$. Thus,
$\pi\tau([v])=\pi(k'x)=kk'[v]$. Suppose
$\delta_2v=\sum_{i=1}^n\nu(v,v_i)v_i$, then
$\delta_1v_i=\nu(v,v_i)v{+}\delta_2v_i$ for $i=1,\cdots,n$, so\\
\hspace*{4mm}$\pi\tau([v])=[\delta_1\delta_2v]=[\delta_1(\sum_{i=1}^n\nu(v,v_i)v_i)]
=[nv+\delta_2(\sum_{i=1}^n\nu(v,v_i)v_i)]=[nv+\delta^2_2v]=n[v]$.\\
Thus, $k'=n/k$.
 \hfill Q.E.D.\vspace{4mm}

\begin{Definition}{\rm
For a connected deformation graph $(G,\nu)$, its volume is defined
as follows. If $G$ has only one vertex, its volume is $0$. If $G$
has more than one vertices, its volume is the number of vertices of
one of its distance components. The characteristic number of
$(G,\nu)$ is $\chi(G,\nu)=|{\rm det}A|$, the absolute value of the
determinant of any of its representation matrix $A$. For a connected
deformable chain graph $(|G|,\nu)$ with rank $>0$, its
characteristic number is
$$\chi(|G|,\nu)=\frac{\Pi_{n=-\infty}^{+\infty}|H_{2n}(|G|)|}
{\Pi_{n=-\infty}^{+\infty}|H_{2n+1}(|G|)|}=
\frac{\Pi_{n=-\infty}^{+\infty}|H^{2n+1}(|G|)|}{\Pi_{n=-\infty}^{+\infty}|H^{2n}(|G|)|},$$
where $|H|$ denotes the cardinality of the finite group $H$.
}\end{Definition}\vspace{4mm}

\begin{Theorem}{\rm For a connected deformation graph
$(G,\nu)$ with volume $n$ and $r(G,\nu)=r$, \\
\hspace*{65mm}$\chi(G,\nu)^2=r^n.$\\ Specifically, when $n$ is odd,
$r(G,\nu)$ is a square number.

For a connected deformable chain graph $(|G|,\nu)$ with rank
$r>0$,\vspace{-2mm}
$$\chi(|G|,\nu)=\chi(G,\nu)r^{\mu},$$
where $\mu=\sum_{k=0}^{+\infty}k(\mu_{2(k{+}s){+}1}-\mu_{2(k{+}s)})$
and $\mu_{k}=$ number of vertices with degree $k$ and $s$ satisfies
that $\mu_i=0$ for $i\!<\!2s$ and $\mu_{2s+1}\!\neq\!0$.
 }\end{Theorem}

Proof. The first equality is obtained from the equality $AA^T=rE$ of
a representation matrix by taking determinant. Let $|\cdot|_1$ be a
representation gradation with representation matrix $A$. By
definition, $\chi(|G|_1,\nu)=|H_0(|G|_1)|=|{\rm det}
A|=\chi(G,\nu)$. Then the second equality of the theorem is a
corollary of Theorem 1.10.  \hfill Q.E.D.\vspace{4mm}

\begin{Theorem}{\rm Let $(G_i,\nu_i)$, $i=1,2$ be two connected
deformation graph with rank $r_i$ and volume $n_i$. Then
$r(G_1{\times} G_2,\nu)=r_1+r_2$ and
$\chi(G_1{\times}G_2,\nu)=(r_1{+}r_2)^{n_1n_2}$.}\end{Theorem}

Proof. By definition. \hfill Q.E.D.\vspace{4mm}

Notice that $\chi(G{\times}H,\nu)$ in the above equality is no
longer a square since the volume of the product graph is $2n_1n_2$.\vspace{4mm}

\begin{Theorem}{\rm Let $(G,\nu)$ be a connected
deformation graph with rank $r>0$ and $F$ be a field of
characteristic $p$. If $p=0$ or $p>0$ but is not a divisor of $r$,
then for all deformable chain graphs $(|G|,\nu)$, $H_*(|G|;F)=0$ and
$H^*(|G|;F)=0$ .}\end{Theorem}

Proof. By Theorem 1.11. \hfill Q.E.D.\vspace{4mm}

\begin{Definition}{\rm A finite graph $G$ is called a diamond graph if it
has the following property. If there are three vertices $a,b,c$ of
$G$ such that $[a,b]$ and $[b,c]$ are edges, then there exists one
and only one new vertex $d$ such that $[c,d]$ and $[d,a]$ are edges
and neither of $[a,c]$ and $[b,d]$ is an edge. The subgraph
consisting of such four vertices and four edges is called a diamond
of $G$. We use four vertices $a,b,c,d$ to denote a diamond such that
$[a,b],[b,c],[c,d],[d,a]$ are edges and $[a,c]$ and $[b,d]$ are not
edges.}\end{Definition}\vspace{4mm}

The above definition implies that there is no triangle in a diamond
graph. \vspace{4mm}

\begin{Theorem}{\rm For a connected diamond graph $G$ , all its vertices
have the same number of neighbors which is called the rank of $G$
and is denoted by $r(G)$. }\end{Theorem}

Proof. If $G$ has no edge, then it has only one vertex with rank 0.
If $a$ is a vertex of $G$ that has rank 1, then $G$ has only one
edge $[a,b]$ and two vertices $a$ and $b$. The number of neighbors
of $a$ and $b$ are all 1 and the conclusion holds. Suppose $G$ has a
vertex $a$ with $n>1$ neighbors. Let $b$ be a neighbor of $a$ and
$A$ and $B$ are respectively the set of neighbors of $a$ and $b$.
For any $v\in A-\{b\}$, three vertices $v,a,b$ determine a unique
vertex $w$ such that $v,a,b,w$ form a diamond. This obviously sets
up a 1-1 correspondence between $A-\{b\}$ and $B-\{a\}$. So $A$ and
$B$ have the same cardinality. Since $G$ is connected, all its
vertices have the same number of neighbors. \hfill
Q.E.D.\vspace{4mm}

\begin{Example}{\rm There exist ungradable diamond graphs. For example, let $D_{1}$
be the diamond graph with vertex set
$\{v,v_i,v_{i,j}{=}v_{j,i}\,|\,1{\leqslant}i{<}j{\leqslant} 5\}$ and
the edges $[v,v_i]$, $[v_i,v_{i,j}]$, $[v_{i,j},v_{s,t}]$ if
$\{i,j\}\cap\{s,t\}=\phi$. The distance function $|u|=d(u,v)$ is not
a gradation, for $|v_{i,j}|=|v_{s,t}|=2$ but $[v_{i,j},v_{s,t}]$ is
an edge for $\{i,j\}\cap\{s,t\}=\phi$. By Theorem 1.3, a distance
function of a gradable graph must be a gradation. So $D_{1}$ is not
gradable.}\end{Example} \vspace{4mm}

\begin{Definition}{\rm
Let $G$ be a gradable diamond graph. A signature $\nu$ on $G$ is a
deformable connection on $G$ such that $\nu(a,b)=\pm1$ for all edges
$[a,b]$. }\end{Definition}\vspace{4mm}

An equivalent definition of a signature is that for every diamond,
three of the four edges have the same sign of signature and the
other edge have the other sign of signature.\vspace{4mm}

\begin{Theorem}{\rm Let $G$ be a diamond graph.
For any deformable connection $\nu$ on $G$, there is a unique
associated signature $\overline\nu$ defined by
$\nu(a,b)=|\nu(a,b)|\overline\nu(a,b)$ for all edges $[a,b]$ and
$\overline\nu(a,b)=0$ if $\nu(a,b)=0$. Two deformable connections
$\nu$ and $\nu'$ are equivalent if and only if their associated
signatures $\overline\nu$ and $\overline\nu'$ are equivalent and
$|\nu(a,b)|=|\nu'(a,b)|$ for all $a,b$. Moreover, all signatures on
$G$ if they exist are equivalent. Thus, a gradable diamond graph is
deformable if and only if there is a signature on it.}\end{Theorem}

Proof. We only prove the uniqueness of the equivalent class of
signatures. Other conclusions are trivial. We may suppose the
diamond graph $G$ is connected. Let $\tilde G$ be the
$2$-dimensional CW-complex defined as follows. The $1$-skeleton of
$\tilde G$ is just the graph $G$ with its usual CW-complex
structure. To every diamond we associated a $2$-cell with attaching
map a homeomorphism from $S^1$ to the four edges of the diamond.
Take a maximal tree on the graph $G$ and suppose $E_1,\cdots,E_n$
are the edges that is not in the maximal tree. Then $\pi_1(G)$
(regard $G$ as a CW-complex) is the free group generated by any $n$
loops that successively containing only one edge $E_i$ for
$i=1,\cdots,n$. Since $G$ is connected, every edge is the edge of a
diamond. This implies that $\pi_1(\tilde G)=0$.

Let $\nu,\nu'$ be two signatures. Take a fixed vertex $v$ of $G$ and
for a path $\omega=\{v,v_1,\cdots,v_n,u\}$, define
$e(\omega,u)=\nu(v,v_1)\nu'(v,v_1)\nu(v_1,v_2)\nu'(v_1,v_2)\cdots\nu(v_n,u)\nu'(v_n,u)$
($\nu(a,a){=}\nu'(a,a){=}1$). It is a direct checking that if two
paths $\omega_1$ and $\omega_2$ differ only on a diamond,
$e(\omega_1,u)=e(\omega_2,u)$. Thus $e$ is invariant on homotopic
loops in $\tilde G$. Since $\pi_1(\tilde G)=0$, $e(u,\omega)$ only
depends on the end vertex $u$ and so $e(u)=e(\omega,u)$ is
well-defined. So for any $a,b\in G$, $\nu(a,b)=e(a)e(b)\nu'(a,b)$.
$\nu$ and $\nu'$ are equivalent.
 \hfill Q.E.D.\vspace{4mm}

\begin{Example}{\rm There exist gradable diamond graphs that has no signature. Let
$D_{2}$ be the diamond graph with vertex set
$\{v,v_i,v_{i,j}{=}v_{j,i},u_1,\cdots,u_6\,|\,1{\leqslant}i{<}j{\leqslant}
5\}$ and gradation $|v|=0$, $|v_i|=1$, $|v_{i,j}|=2$, $|u_k|=3$. The
edges are $[v,v_i]$, $[v_i,v_{i,j}]$, and\\
\hspace*{4mm}$[u_1,v_{1,2}]$,$[u_1,v_{2,3}]$,$[u_1,v_{3,4}]$,$[u_1,v_{4,5}]$,$[u_1,v_{5,1}]$,
$[u_2,v_{1,2}]$,$[u_2,v_{2,4}]$,$[u_2,v_{4,5}]$,$[u_2,v_{5,3}]$,$[u_2,v_{3,1}]$,\\
\hspace*{4mm}$[u_3,v_{1,2}]$,$[u_3,v_{2,5}]$,$[u_3,v_{5,3}]$,$[u_3,v_{3,4}]$,$[u_3,v_{4,1}]$,
$[u_4,v_{3,2}]$,$[u_4,v_{2,4}]$,$[u_4,v_{4,1}]$,$[u_4,v_{1,5}]$,$[u_4,v_{5,3}]$,\\
\hspace*{4mm}$[u_5,v_{3,2}]$,$[u_5,v_{2,5}]$,$[u_5,v_{5,4}]$,$[u_5,v_{4,1}]$,$[u_5,v_{1,3}]$,
$[u_6,v_{4,2}]$,$[u_6,v_{2,5}]$,$[u_6,v_{5,1}]$,$[u_6,v_{1,3}]$,$[u_6,v_{3,4}]$.\\
The volume of $D_{2}$ is $11$ but the rank is $5$ and not a square
number. So by Theorem 1.11, $D_{2}$ has no deformable connection on
it. }\end{Example}\vspace{4mm}

\begin{Definition}{\rm
A gradable diamond graph $G$ that has a signature is called
admissible. For a gradation $|\cdot|$ on $G$, the deformable chain
graph $(|G|,\nu)$ with $\nu$ any of its signature is called a GAD
(graded, admissible, diamond) graph. Since there is only one
equivalent class of signatures, $(|G|,\nu)$ is often simply denoted
by $|G|$. }\end{Definition}\vspace{4mm}

It is obvious that the rank of a connected GAD graph $|G|$ as a
deformable chain graph equals the rank of the diamond graph as
defined in Theorem 1.14.\vspace{4mm}

\begin{Example}{\rm The (co)homology of the distance graph of an admissible
diamond graph may not be trivial. Let $(C,d)$ be defined as follows.
The vertex set is
$\{v,v_i,v_{i,j}{=}{-}v_{j,i},1{\leqslant}i{<}j{\leqslant}
4,$\\$e_1,e_2,e_3\}$, $|v|=0$, $|v_i|=1$, $|v_{i,j}|=2$, $|e_k|=3$,
$dv=0$, $dv_i=v$, $dv_{i,j}=v_i-v_j$,
$de_1=v_{1,2}{+}v_{2,3}{+}v_{3,4}{+}v_{4,1}$,
$de_2=v_{1,3}{+}v_{3,4}{+}v_{4,2}{+}v_{2,1}$,
$de_3=v_{1,4}{+}v_{4,2}{+}v_{2,3}{+}v_{3,1}$. It has the positive
distance gradation of a diamond graph. $H_2(C)=\langle 2\rangle$
with generator class represented by $v_{2,3}{+}v_{3,4}{+}v_{4,2}$;
$H_s(C)=0$, otherwise.}\end{Example}\vspace{4mm}

\section{Diamond root systems}\vspace{4mm}

\hspace{6mm}Recall that a Lie algebra $\frak G$ (over $\Bbb Z$) is
an Abellian group with Lie bracket $[\,,]\colon\frak G\otimes\frak
G\to \frak G$ that satisfies the following properties.

1) $[x,x]=0$ for all $x\in\frak G$.

2) $[[x,y],z]{+}[[y,z],x]{+}[[z,x],y]=0$ for all $x,y,z\in\frak
G$.\vspace{4mm}

In this paper, we only study free Lie algebras, i.e., Lie algebras
that are free Abellian groups generated by a finite set $S$ which is
called the base of the Lie algebra.

\begin{Definition}{\rm
For a Lie algebra $\frak G$ with base $S$, its associated chain
graph $(|\Lambda(S)|,\nu)$ with respect to $S$ is defined as
follows.  Suppose $S=\{e_1,\cdots,e_n\}$ and
$[e_i,e_j]=\sum_{k=1}^nr_{i,j}^ke_k$ with $r_{i,j}^k\in\Bbb Z$ for
$i<j$. Let $\Lambda(\frak G)$ be the exterior algebra generated by
$\frak G$. Then $\Lambda(S)$ is the base of $\Lambda(\frak G)$
consisting of $1$ and elements of the form $e_{i_1}e_{i_2}\cdots
e_{i_s}$ with
$1\,{\leqslant}\,i_1\,{<}\,i_2\,{<}\cdots{<}i_s\,{\leqslant}\,n$
with gradation $|1|=0$ and $|e_{i_1}e_{i_2}\cdots e_{i_s}|=s$. There
are two dual derivatives $d$ and $\delta$ on $\Lambda(\frak G)$.
$(\Lambda(\frak G),\delta)$ is a DGA such that $\delta
e_i=\sum_{s<t}r_{s,t}^i e_se_t$ for $i=1,\cdots,n$ and
$d(e_{i_1}e_{i_2}\cdots
e_{i_s})=\sum_{u<v}(-1)^{v-u}[e_{i_u},e_{i_v}]e_{i_1}\cdots \hat
e_{i_u}\cdots \hat e_{i_v}\cdots e_{i_s}$. $d$ and $\delta$
naturally induce the same connection $\nu$ on $|\Lambda(S)|$. Thus,
$(|\Lambda(S)|,\nu)$ is the chain graph with respectively associated
chain complex and cochain complex $(\Lambda(\frak G),d)$ and
$(\Lambda(\frak G),\delta)$. $\Lambda(S)$ ($|\Lambda(S)|$) is called
the (graded) base graph of $\frak G$ and $\nu$ is called the
associated connection of $\frak G$ with respect to base $S$.
}\end{Definition}\vspace{4mm}

\begin{Definition}{\rm
A diamond Lie algebra $\frak G$ is a Lie algebra with base $S$ that
satisfies the following property. Such a base $S$ is called a
diamond root system.

1) For $\alpha\in S$, $-\alpha\not\in S$.

2) For two different $\alpha,\beta\in S$, either $[\alpha,\beta]=0$
or $\pm[\alpha,\beta]\in S{-}\{\alpha,\beta\}$.

3) For three different $\alpha,\beta,\gamma\in S$ such that
$\pm[\alpha,\beta],\pm[\beta,\gamma]\in S$, $[\alpha,\gamma]=0$ and
$\pm[\alpha,\beta,\gamma]=\pm[[\alpha,\beta],\gamma]=\pm[\alpha,[\beta,\gamma]]\in
S$ and $[\alpha,\beta]\neq\pm[\beta,\gamma]$. Such three different
$\alpha,\beta,\gamma$ are called adjacent.

4) For four different $\xi,\eta,\sigma,\tau\in S$ such that
$[\xi,\eta]=[\sigma,\tau]\in S$, there are adjacent
$\alpha,\beta,\gamma\in S$ such that $\xi=\alpha$,
$\eta=\pm[\beta,\gamma],\,\sigma=\pm[\alpha,\beta]$, $\tau=\gamma$
(permute $\xi,\eta,\sigma,\tau$ if necessary) and there are no
adjacent $\alpha',\beta',\gamma'\in S$ such that $\eta=\alpha'$,
$\xi=\pm[\beta',\gamma']$. }\end{Definition}\vspace{4mm}

\begin{Theorem}{\rm For a Lie algebra $\frak G$ with diamond root system $S$, its base
graph $\Lambda(S)$ is a diamond graph and the associated chain graph
$(|\Lambda(S)|,v)$ is a GAD graph.}\end{Theorem}

Proof. We will prove that for a diamond root system $S$,
$\Lambda(S)$ is a  diamond graph with only the following six types
of diamond
($[\alpha,\beta,\gamma]=[[\alpha,\beta],\gamma]=[\alpha,[\beta,\gamma]]$
and an arrow represents an edge from the vertex with higher degree
to the vertex with lower degree)
\begin{eqnarray*}
(1)&\left.\begin{array}{c}\xi\eta\sigma\tau x\\
\swarrow\quad\quad\searrow\\
{[}\xi,\eta{]}\sigma\tau x\quad\quad\xi\eta{[}\sigma,\tau{]} x\\
\searrow\quad\quad\swarrow\\
{[}\xi,\eta{]}{[}\sigma,\tau{]} x
\end{array}\right.\\ \\
(2)&\left.\begin{array}{c}\alpha\beta\gamma x\\
\swarrow\quad\quad\searrow\\
\alpha[\beta,\gamma] x\quad\quad\gamma[\alpha,\beta] x\\
\searrow\quad\quad\swarrow\\
{[}\alpha,\beta,\gamma{]} x\end{array}\right.\\ \\
(3)&\left.\begin{array}{c}\alpha\beta\gamma{[}\alpha,\beta{]} x\\
\swarrow\quad\quad\quad\searrow\\
\alpha{[}\alpha,\beta{]}{[}\beta,\gamma{]} x\quad\quad\alpha\beta{[}\alpha,\beta,\gamma{]} x\\
\searrow\quad\quad\quad\swarrow\\
{[}\alpha,\beta{]}{[}\alpha,\beta,\gamma{]} x\\
\end{array}\right.\\ \\
(4)&\left.\begin{array}{c}\alpha\beta\gamma{[}\alpha,\beta{]}{[}\beta,\gamma{]} x\\
\swarrow\quad\quad\quad\quad\searrow\\
\alpha\beta{[}\beta,\gamma{]}{[}\alpha,\beta,\gamma{]} x\quad\quad
\beta\gamma{[}\alpha,\beta{]}{[}\alpha,\beta,\gamma{]} x\\
\searrow\quad\quad\quad\quad\swarrow\\
{[}\alpha,\beta{]}{[}\beta,\gamma{]}{[}\alpha,\beta,\gamma{]} x\\
\end{array}\right.\\ \\
(5)&\left.\begin{array}{c}\alpha\beta[\beta,\gamma]x\quad\quad\beta\gamma[\alpha,\beta]x\\
\downarrow\quad\quad\quad\swarrow\hspace{-4mm}\searrow\quad\quad\quad\downarrow\\
{[}\alpha,\beta{]}{[}\beta,\gamma{]}x\quad\quad\beta{[}\alpha,\beta,\gamma{]}x\\
\end{array}\right.\\ \\
(6)&\left.\begin{array}{c}\alpha\gamma[\alpha,\beta][\beta,\gamma]x\quad\quad\alpha\beta\gamma[\alpha,\beta,\gamma]x\\
\downarrow\quad\quad\quad\swarrow\hspace{-4mm}\searrow\quad\quad\quad\downarrow\\
\alpha[\beta,\gamma][\alpha,\beta,\gamma{]}x\quad\quad\gamma[\alpha,\beta][\alpha,\beta,\gamma]x\\
\end{array}\right.
\end{eqnarray*}
where $x$ is a product of elements of $S$ with no factors appearing
in the front.

The uniqueness of the above six types of diamonds is from the
definition of the diamond system.

In the following discussion, sign is neglected. That is, $x=y$
implies $x=\pm y$. Thus, we may suppose the root system satisfies
that $[\alpha,\beta]=\alpha{+}\beta$ for all roots
$\alpha,\beta,\alpha{+}\beta$.

For three $u,v,w\in \Lambda(S)$ such that $u$ and $v$ are neighbors
and $v$ and $w$ are neighbors, they are one of the following, where
$\xi,\eta,\sigma,\tau\in S$ and $x,y\in \Lambda(S)$.

1) $u=[\xi,\eta]x$, $v=\xi\eta x=[\sigma,\tau]y$, $w=\sigma\tau y$;

2) $u=[\xi,\eta]x$, $v=\xi\eta x=\sigma\tau y$, $w=[\sigma,\tau]y$;

3) $u=\xi\eta x$, $v=[\xi,\eta]x=[\sigma,\tau]y$, $w=\sigma\tau y$.

If the six elements $\xi,\eta,\sigma,\tau,[\xi,\eta],[\sigma,\tau]$
of $S$ are different and so $x$ and $y$ have none of the six
factors, then all the above cases are three vertices of the diamond
of type (1). If there are repetitions in the six elements, then from
the symmetry of the pairs $\xi$, $\eta$ and $\sigma$, $\tau$, we
need only prove the following four cases. (a)
$[\xi,\eta]=[\sigma,\tau]$; (b) $\eta=[\sigma,\tau]$; (c)
$\tau=[\xi,\eta]$; (d) $\xi=\sigma$.

(a) In this case, there are three $\alpha,\beta,\gamma\in S$ such
that $\xi=\alpha$, $\eta=[\beta,\gamma],\,\sigma=[\alpha,\beta]$,
$\tau=\gamma$. If $u,v,w$ are case 1), then
$v=\xi\eta[\sigma,\tau]z$ and so $u=[\xi,\eta][\sigma,\tau]z=0$.
Impossible. If $u,v,w$ are case 2), then
$u=\gamma[\alpha,\beta][\alpha,\beta,\gamma]z$,
$v=\alpha\gamma[\alpha,\beta][\beta,\gamma]z$,
$w=\alpha[\beta,\gamma][\alpha,\beta,\gamma]z$. If $z$ has no factor
$\beta$, then $u,v,w$ are three vertices of a diamond of type $(6)$.
If $z$ has factor $\beta$, then $u,v,w$ are three vertices of a
diamond of type $(4)$. If $u,v,w$ are case 3), then
$u=\alpha[\alpha,\beta,\gamma]z$, $v=[\alpha,\beta,\gamma]z$,
$w=\gamma[\alpha,\beta]z$. If $z$ has no factor $\beta$, then
$u,v,w$ are three vertices of a diamond of type $(2)$. If $z$ has
factor $\beta$, $u,v,w$ are three vertices of a diamond of type
$(5)$.

(b) In this case, $\eta=[\sigma,\tau]$. By Jacobi identity,
$[\xi,[\sigma,\tau]]{+}[\sigma,[\tau,\xi]]{+}[\tau,[\xi,\sigma]]=0$,
only one of $[\xi,\sigma]$ and $[\xi,\tau]$ is not zero. We may
suppose $[\xi,\sigma]\neq 0$ and $[\xi,\tau]=0$. Thus,
$\xi,\sigma,\tau$ are adjacent. If $u,v,w$ are case 1), then
$u=[\xi,\sigma,\tau]z$, $v=\xi[\sigma,\tau]z$, $w=\xi\sigma\tau z$.
If $z$ has no factor $[\xi,\sigma]$, then $u,v,w$ are three vertices
of a diamond of type $(2)$. If $z$ has factor $[\xi,\sigma]$,
$u,v,w$ are three vertices of a diamond of type $(3)$. If $u,v,w$
are case 2), then $v=\xi\eta\sigma\tau z$, $w=
\xi\eta[\sigma,\tau]z=0$. Impossible. If $u,v,w$ are case 3), then
$v=[\xi,\eta][\sigma,\tau]z$, $u=\xi\eta[\sigma,\tau]z=0$.
Impossible.

(c) In this case, $\sigma=[\xi,\eta]$. By Jacobi identity,
$[[\xi,\eta],\tau]]{+}[[\eta,\tau],\xi]{+}[[\tau,\xi],\eta]=0$, only
one of $[\xi,\tau]$ and $[\eta,\tau]$ is not zero. We may suppose
$[\xi,\tau]=0$ and $[\eta,\tau]\neq 0$. Thus, $\xi,\eta,\tau$ are
adjacent. If $u,v,w$ are case 1), then
$u=[\xi,\eta][\xi,\eta,\tau]z$, $v=\xi\eta[\xi,\eta,\tau]z$,
$w=\xi\eta\tau[\xi,\eta]z$. If $z$ has no factor $[\eta,\tau]$, then
$u,v,w$ are three vertices of a diamond of type $(3)$. If $z$ has
factor $[\eta,\tau]$, $u,v,w$ are three vertices of a diamond of
type $(4)$. If $u,v,w$ are case 2), then $v=\xi\eta\sigma\tau z$,
$u=[\xi,\eta]\sigma\tau z=0$. Impossible. If $u,v,w$ are case 3),
then $v=[\xi,\eta][\sigma,\tau]z$, $w=[\xi,\eta]\sigma\tau z=0$.
Impossible.

(d) In this case, $\xi=\sigma$, $\eta\neq\tau$ and $\eta,\xi,\tau$
are adjacent. If $u,v,w$ are case 1), then $v=\xi\eta[\xi,\tau]z$
and so $w=\xi\eta\xi\tau z=0$. Impossible. If $u,v,w$ are case 2),
then $u=[\xi,\eta]\tau z$, $v=\xi\eta\tau z$, $w=\eta[\xi,\tau]z$.
If $z$ has no factor $[\eta,\xi,\tau]$, then $u,v,w$ are three
vertices of a diamond of type $(2)$. If $z$ has factor
$[\eta,\xi,\tau]$, then $u,v,w$ are three vertices of a diamond of
type $(6)$. If $u,v,w$ are case 3), then $u=\xi\eta[\xi,\tau]z$,
$v=[\xi,\eta][\xi,\tau]z$, $w=\xi\tau[\eta,\tau]z$. If $z$ has no
factor $[\eta,\xi,\tau]$, then $u,v,w$ are three vertices of a
diamond of type $(5)$. If $z$ has factor $[\eta,\xi,\tau]$, $u,v,w$
are three vertices of a diamond of type $(4)$.

Overall, the types of diamond corresponding to case 1),2),3) and
(a),(b),(c),(d) are as in the following table.
$$\begin{array}
{|c|c|c|c|} \hline &1)&2)&3)\\
\hline (a)&{\rm impossible}&(4)\,{\rm or}\,(6)&(2)\,{\rm or}\,(5)\\
\hline (b)&(2)\,{\rm or}\,(3)&{\rm impossible}&{\rm impossible}\\
\hline (c)&(3)\,{\rm or}\,(4)&{\rm impossible}&{\rm impossible}\\
\hline (d)&{\rm impossible}&(2)\,{\rm or}\,(6)&(4)\,{\rm or}\,(5)\\
\hline
\end{array}$$
 \hfill Q.E.D.\vspace{4mm}

\begin{Theorem}{\rm Let $\frak G$ be the Lie algebra (over $\Bbb Z$)
generated by the positive root system of a simple Lie algebra over
$\Bbb C$. Except the case when the Dynkin graph is $C_n$ ($n>2$) or
$F_4$, $\frak G$ is a diamond Lie algebra.}\end{Theorem}

Proof. Direct checkings. \hfill Q.E.D.\vspace{4mm}

The chain graph $(|\Lambda(\frak G)|,\nu)$ when the Dynkyn graph is
$C_n$ or $F_4$ is even not a deformable chain graph.\vspace{4mm}

For a diamond root system $S$, give $S$ an order and define $w(S)$
to be the Abellian group generated by $S$ modular zero relations
$[\alpha,\beta]=\alpha{+}\beta$ for all $\alpha<\beta$ and
$[\alpha,\beta]\neq 0$. A tedious proof shows that if $w(S)$ is a
free Abellian group generated by $T$ such that $T$ is a connected
graph ($[a,b]$ is an edge of $T$ if and only if $a+b\in S$), then
$S$ is one of those in Theorem 2.2.\vspace{4mm}

\begin{Definition}{\rm Let $\frak G$  be a diamond Lie algebra with
base $S$. Let $\omega(S)$ denote the set of connected components of
the base graph $\Lambda(S)$. Then there is a graph connected
component decomposition
$\Lambda(S)=\sqcup_{\alpha\in\omega(S)}\Lambda(\alpha)$ and
corresponding chain complex and cochain
complex decompositions\\
\hspace*{26mm}$(\Lambda(\frak
G),d)=\oplus_{\alpha\in\omega(S)}(\Lambda(\alpha),d),\quad\quad
(\Lambda(\frak
G),\delta)=\oplus_{\alpha\in\omega(S)}(\Lambda(\alpha),\delta).$

Denote $H_{\alpha,*}(\frak G)=H_{*}(\Lambda(\alpha),d)$ and
$H^{\alpha,*}(\frak G)=H^{*}(\Lambda(\alpha),\delta)$, then there is
a direct sum decomposition\\
\hspace*{31mm}$H_*(\frak
G)=\oplus_{\alpha\in\omega(S)}H_{\alpha,*}(\frak G),\quad\quad
H^*(\frak G)=\oplus_{\alpha\in\omega(S)}H^{\alpha,*}(\frak G).$
}\end{Definition}\vspace{4mm}

\begin{Theorem}{\rm Let $\frak G$  be a diamond Lie algebra with
base $S$ and for $\alpha\in\omega(S)$, $r(\alpha)$ denote the rank
of the connected diamond graph $\Lambda(\alpha)$. Then the free part
of the (co)homology group of $\frak G$ is the free Abellian group
generated by all $\alpha\in\omega(S)$ with $r(\alpha)=0$. For a
prime $p$, if $r(\alpha)$ is not divisible by $p$, then
$H_{\alpha,*}(\frak G)$ and $H^{\alpha,*}(\frak G)$ have no
$p$-torsion part. }\end{Theorem}

Proof. By Theorem 1.13. \hfill Q.E.D.\vspace{4mm}

\section{Weight (co)chain subcomplexes of ${\frak A}_{n+1}$}\vspace{4mm}

\hspace*{6mm}Let $\frak A_{n+1}$ be the Lie algebra generated by the
positive root system with Dynkin graph $A_{n+1}$. We will determine
all the connected components of $\Lambda(\frak A_{n+1})$ and compute
their rank. Precisely, let $\frak A_{n+1}$ be the Lie algebra of all
$(n{+}1){\times}(n{+}1)$ upper triangular matrices with integer
entries and diagonal zero. For simplicity, we denote the entries of
a $(n{+}1){\times}(n{+}1)$ matrix by $a_{i,j}$ with
$i,j=0,1,\cdots,n$. Let $\{e_{i,j}\,|\,0\leqslant i<j\leqslant n\}$
denote the matrix with all entries $0$ except $a_{i,j}=1$, then the
Lie bracket of the Lie algebra is defined by
$$[e_{i,j},e_{k,l}]=\left\{\begin{array}{cl}e_{i,l}&{\rm if}\,\,j=k\\
-e_{j,k}&{\rm if}\,\,i=l\\0\,&{\rm otherwise}\end{array}\right..$$
We denote the associated chain and cochain complexes $(\Lambda(\frak
A_{n+1}),d)$ and $(\Lambda(\frak A_{n+1}),\delta)$ respectively by
$(L_n,d)$ and $(R_n,\delta)$.\vspace{4mm}

\begin{Definition}{\rm For $n=1,2,\cdots$, let $G_n$ be
the set of all triangular matrix $[a_{i,j}]$ with entries
$a_{i,j}=0$ or $1$
$$[a_{i,j}]=\left[\begin{array}{llllc}
a_{0,1}&a_{0,2}&a_{0,3}&\cdots&a_{0,n}\\
&a_{1,2}&a_{1,3}&\cdots&a_{1,n}\\
&&a_{2,3}&\cdots&a_{2,n}\\
&&&\cdots&\cdots\\
&&&&a_{n-1,n}
\end{array}\right].$$
For a triangular matrix $[a_{i,j}]\in G_n$, its weight
$(i_0,i_1,\cdots,i_n)$ is defined to be
$i_s=\sum_{k=0}^{s-1}(1-a_{k,s})+\sum_{k=s+1}^n a_{s,k}$,
$s=0,1,\cdots,n$. $G(i_0,i_1,\cdots,i_n)$ is defined to be the
subset of $G_n$ of all triangular matrices with weight
$(i_0,i_1,\cdots,i_n)$. Specifically, we define $G_0=\{[\,]\}$
having only one `vacuum' triangular matrix $[\,]$ with weight $(0)$.

We correspond a triangular matrix $[a_{i,j}]$ to a monomial
$e_{0,1}^{a_{0,1}}e_{1,2}^{a_{1,2}}e_{0,2}^{a_{0,2}}\cdots
e_{n-1,n}^{a_{n-1,n}}\cdots e_{1,n}^{a_{1,n}}e_{0,n}^{a_{0,n}}$
$(e_{i,j}^1=e_{i,j}\, ,e_{i,j}^0=1)$ in $\Lambda(\frak A_{n+1})$.
With this correspondence, $G_n=\Lambda(\frak A_{n+1})$. We denote by
$G(i_0,i_1,\cdots,i_n)$ the set of all triangular matrices with
weight $(i_0,i_1,\cdots,i_n)$ and call it the weight subgraph of
$G_n$. We denote the chain complex $(G(i_0,i_1,\cdots,i_n),d)$ and
the cochain complex $(G(i_0,i_1,\cdots,i_n),\delta)$ respectively by
$(L(i_0,i_1,\cdots,i_n),d)$ and $(R(i_0,i_1,\cdots,i_n),\delta)$ and
call them the weight chain and cochain
subcomplexes.}\end{Definition}\vspace{4mm}

{\bf Remark.} We use a triangular matrix $[a_{i,j}]$ to denote both
the vertex of $G_n$ and the product $\prod_{i,j}e_{i,j}^{a_{i,j}}$
in $\Lambda(\frak A_{n+1})$ as defined in Definition 3.1. Since
$R_n$ is a DGA, it is sometimes convenient to discuss the problem on
$R_n$ but not on $G_n$ or $L_n$. So we denote the algebra generator
of $R_n$ by $x_{i,j}$ to distinguish them from $e_{i,j}$ in $G_n$
and $L_n$.\vspace{4mm}

\begin{Theorem}{\rm The weight subgraphs satisfy the following.
\begin{enumerate}
\item Let $S_{n+1}$ be the
group of permutations on $\{0,1,\cdots,n\}$. For any $\sigma\in
S_{n+1}$, there is a graph isomorphism $g(\sigma):G_n\to G_n$ such
that the restriction map $g(\sigma)|_{G(i_0,i_1,\cdots,i_n)}$ is a
graph isomorphism from $G(i_0,i_1,\cdots,i_n)$ to
$G(i_{\sigma(0)},i_{\sigma(1)},\cdots,i_{\sigma(n)})$ which is
generally not a GAD graph isomorphism.
\item  For any weight $(i_0,i_1,\cdots,i_{n-1},i_n)$, there is a transpose
GAD graph isomorphism from $G(i_0,i_1,\cdots,i_{n-1},i_n)$ to
$G(n{-}i_n,n{-}i_{n-1},\cdots,n{-}i_1,n{-}i_0)$.
\item  For any weight $(i_0,i_1,\cdots,i_{n-1},i_n)$, there is a rotation GAD graph
isomorphism from $G(i_0,i_1,\cdots,i_{n-1},i_n)$ to
$G(i_n,i_0,i_1,\cdots,i_{n-1})$.
\item  For any weight $(i_0,i_1,\cdots,i_{n-1},i_n)$,
there is a duality GAD graph isomorphism from
$G(i_0,i_1,\cdots,i_{n-1},i_n)$ to
$G(n{-}i_0,n{-}i_1,\cdots,n{-}i_{n-1},n{-}i_n)$.\end{enumerate}
}\end{Theorem}

Proof. 1. For $k=1,\cdots,n$, let $\sigma_k=(k{-}1,k)\in S_{n+1}$
($\sigma_k(k{-}1)=k$, $\sigma_k(k)=k{-}1$, $\sigma_k(i)=i$ for
$i\neq k{-}1,k$). Define $g(\sigma_k)$ as follows. For any
triangular matrices $[a_{i,j}]\in G_n$,
\begin{eqnarray*}&g(\sigma_k)( \left[\begin{array}{cccccc}
\cdots&a_{0,k-1}&a_{0,k}&\cdots&\cdots&\cdots\\
\cdots&\cdots&\cdots&\cdots&\cdots&\cdots\\
&a_{k-2,k-1}&a_{k-2,k}&\cdots&\cdots&\cdots\\
&&a_{k-1,k}&a_{k-1,k+1}&\cdots&a_{k-1,n}\\
&&&a_{k,k+1}&\cdots&a_{k,n}\\
&&&&\cdots&\cdots
\end{array}\right])\\
=&\quad\quad\quad\left[\begin{array}{cccccc}
\cdots&a_{0,k}&a_{0,k-1}&\cdots&\cdots&\cdots\\
\cdots&\cdots&\cdots&\cdots&\cdots&\cdots\\
&a_{k-2,k}&a_{k-2,k-1}&\cdots&\cdots&\cdots\\
&&1-a_{k-1,k}&a_{k,k+1}&\cdots&a_{k,n}\\
&&&a_{k-1,k+1}&\cdots&a_{k-1,n}\\
&&&&\cdots&\cdots
\end{array}\right],\end{eqnarray*}
where the omitted part remains unchanged. It is a direct checking
that if a triangular matrix $a$ has weight
$(\cdots,i_{k-1},i_k,\cdots)$, then $g(\sigma_k)(a)$ has weight
$(\cdots,i_k,i_{k-1},\cdots)$. Since $\sigma_1,\cdots,\sigma_n$
generate $S_{n+1}$, $g(\sigma)$ is defined for all $\sigma\in
S_{n+1}$.

2. Define DGA isomorphism $\varrho\colon R_n\to R_n$ by
$\lambda(x_{i,j})=x_{n-j,n-i}$ . Then for any triangular matrix
$[a_{i,j}]\in G_n$,
\begin{eqnarray*}&\varrho(\left[\begin{array}{rrrrr}
a_{0,1}&a_{0,2}&\cdots&a_{0,n-1}&a_{0,n}\\
&a_{1,2}&\cdots&a_{1,n-1}&a_{1,n}\\
&&\cdots&\cdots&\cdots\\
&&&a_{n-2,n-1}&a_{n-2,n}\\
&&&&a_{n-1,n}
\end{array}\right])\\
=&\pm\left[\begin{array}{lllll}
a_{n-1,n}&a_{n-2,n}&\cdots&a_{1,n}&a_{0,n}\\
&a_{n-2,n-1}&\cdots&a_{1,n-1}&a_{0,n-1}\\
&&\cdots&\cdots&\cdots\\
&&&a_{1,2}&a_{0,2}\\
&&&&a_{0,1}
\end{array}\right]\end{eqnarray*}
and $\varrho$ induces a GAD graph isomorphism that sends
$G(i_0,\cdots,i_n)$ to $G(n{-}i_n,\cdots,n{-}i_0)$.

3. Notice that $R_n$ has two DGA-module structures over $R_{n-1}$.
The first module structure is induced by $R_{n-1}$  as a subalgebra
and in this sense, $R_n$ is freely generated by the set
$\{1\}\cup\{x_{i_0,n}\cdots x_{i_s,n},0\leqslant i_0<\cdots< i_s<
n\}$. The second is induced by the monomorphism $j\colon R_{n-1}\to
R_n$ defined by $j(x_{i,j})=x_{i+1,j+1}$ and in this sense, $R_n$ is
freely generated by the set $\{1\}\cup\{x_{0,i_0}\cdots
x_{0,i_s},0<i_0<\cdots< i_s\leqslant n\}$. Define $\varpi\colon
R_n\to R_n$ as follows. For any $a\in R_{n-1}$ and $0\leqslant
i_1<\cdots< i_s<n$, $\varpi(ax_{i_0,n}\cdots
x_{i_s,n})=(-1)^{i_0+\cdots i_s}j(a)x_{0,1}\cdots\hat
x_{0,i_0+1}\cdots\hat x_{0,i_s+1}\cdots x_{0,n}$ (the hat represents
cancelling the factor in the product $x_{0,1}x_{0,2}\cdots x_{0,n}$
and the term is abbreviated to $\cdots\hat x_{0,i_0+1}\cdots\hat
x_{0,i_s+1}\cdots$ in the following formulas). Then,
\begin{eqnarray*}&&\varpi\left(\delta(ax_{i_s,n}\cdots x_{i_0,n})\right)\\
&=&\varpi\left((\delta a)(x_{i_s,n}\cdots x_{i_0,n})\right)\\
&&+\sum_{k,j}(-1)^{|a|+l-k}\varpi(ax_{i_k,j}x_{i_s,n}\cdots
x_{i_{l-1},n}x_{j,n}x_{i_l,n}\cdots
\hat x_{i_k,n}\cdots x_{i_0,n})\\
&=&(-1)^{i_0+\cdots+i_s}(j(\delta a))(\cdots\hat x_{0,i_0+1}\cdots\hat x_{0,i_s+1}\cdots)\\
&&+\sum_{k,j}(-1)^{|a|+l-k+i_0+\cdots+i_s+j-i_k}(j(ax_{i_k,j}))(\cdots\hat
x_{0,i_0+1}\cdots\hat x_{0,j+1}\cdots
x_{i_k,n}\cdots\hat x_{0,i_s+1}\cdots)\\
&=&(-1)^{i_0+\cdots+i_s}(\delta(ja))(\cdots\hat x_{0,i_0+1}\cdots\hat x_{0,i_s+1}\cdots)\\
&&+\sum_{k,j}(-1)^{|a|+l-k+i_0+\cdots+i_s+j-i_k}(ja)(x_{i_k+1,j+1}\cdots\hat
x_{0,i_0+1}\cdots\hat x_{0,j+1}\cdots
x_{i_k,n}\cdots\hat x_{0,i_s+1}\cdots)\\
&=&(-1)^{i_0+\cdots+i_s}(\delta(ja))(\cdots\hat x_{0,i_0+1}\cdots\hat x_{0,i_s+1}\cdots)\\
&&+\sum_{k,j}(-1)^{|a|+i_0+\cdots+i_s}(ja)(\delta(\cdots\hat
x_{0,i_0+1}\cdots\hat x_{0,i_s+1}\cdots))\\
&=&\delta\left(\varpi(ax_{i_s,n}\cdots x_{i_0,n})\right).
\end{eqnarray*}
Thus, $\varpi$ is a DGA-module isomorphism satisfying that for any
triangular matrix $[a_{i,j}]\in G_n$,
\begin{eqnarray*}&\varpi(\left[\begin{array}{rrrrr}
a_{0,1}&a_{0,2}&\cdots&a_{0,n-1}&a_{0,n}\\
&a_{1,2}&\cdots&a_{1,n-1}&a_{1,n}\\
&&\cdots&\cdots&\cdots\\
&&&a_{n-2,n-1}&a_{n-2,n}\\
&&&&a_{n-1,n}
\end{array}\right])\\
=&\pm\left[\begin{array}{ccccc}
1-a_{0,n}&1-a_{1,n}&1-a_{2,n}&\cdots&1-a_{n-1,n}\\
&a_{0,1}&a_{0,2}&\cdots&\quad a_{0,n-1}\\
&&a_{1,2}&\cdots&\quad a_{1,n-1}\\
&&&\cdots&\cdots\\
&&&&a_{n-2,n-1}
\end{array}\right]\end{eqnarray*}
and induces a GAD graph isomorphism that sends $G(i_0,\cdots,i_n)$
to $G(i_n,i_0,\cdots,i_{n-1})$.

4. Define duality isomorphism $\vartheta\colon R_n\to L_n$ as
follows. For any $[a_{i,j}]\in G_n$,
\begin{eqnarray*}&\vartheta(\left[\begin{array}{rrrrr}
a_{0,1}&a_{0,2}&\cdots&a_{0,n-1}&a_{0,n}\\
&a_{1,2}&\cdots&a_{1,n-1}&a_{1,n}\\
&&\cdots&\cdots&\cdots\\
&&&a_{n-2,n-1}&a_{n-2,n}\\
&&&&a_{n-1,n}
\end{array}\right])\\
=&(-1)^{\tau}\left[\begin{array}{lllll}
1-a_{0,1}&1-a_{0,2}&\cdots&1-a_{0,n-1}&1-a_{0,n}\\
&1-a_{1,2}&\cdots&1-a_{1,n-1}&1-a_{1,n}\\
&&\cdots&\cdots&\cdots\\
&&&1-a_{n-2,n-1}&1-a_{n-2,n}\\
&&&&1-a_{n-1,n}
\end{array}\right].\end{eqnarray*}
where $\tau=\sum_{i,j}a_{i,j}[(1+2+\cdots+(j-1)+(j-i-1)]$ . It is
obvious that $\vartheta \delta=d\vartheta$ and so $\vartheta$
induces a GAD graph isomorphism that sends $G(i_0,\cdots,i_n)$ to
$G(n{-}i_0,\cdots,n{-}i_n)$. \hfill Q.E.D.\vspace{4mm}

{\bf Remark} The conclusion 1. and 4. of Theorem 3.1 can be
generalized to all the semi-simple Lie algebras over $\Bbb C$.
Precisely, the Weyl group of a positive system acts on the
associated chain graph $\Lambda(\frak G)$ of the Lie algebra $\frak
G$ generated by the positive system and there is Poncar$\acute{\rm
e}$ duality on it. This is a generalization of Kostant theorem from
the complex number case to the ring of integers case. But other
isomorphisms can not be naturally generalized to even diamond Lie
algebras.\vspace{4mm}

\begin{Definition}{\rm For $n\geqslant 0$, $\omega_n$ is the set
of all $(n{+}1)$-tuples $(i_0,\cdots,i_n)$ such that for all
$0\!\leqslant\!k_0\!<\!\cdots\!<\!k_s\!\leqslant\!n$,
$i_{k_0}\!+\!i_{k_1}\!+\!\cdots \!+\!i_{k_s}\geqslant
0\!+\!1\!+\!\cdots\! +\!s$ and
$i_0\!+\!i_1\!+\!\cdots\!+\!i_n=0\!+\!1+\!\cdots\!+\!n$. For $s>0$,
$(i_0,\cdots,i_{s})\in\omega_s$ is reducible if there is $0\leqslant
m<s$ and $0\leqslant u_0<\cdots<u_m\leqslant s$ such that
$i_{u_0}+\cdots+i_{u_m}=0+1+\cdots +m$. Equivalently,
$(i_0,\cdots,i_{s})\in\omega_s$ is reducible if it is a permutation
of $(j_0,\cdots,j_m,k_0{+}m{+}1,\cdots,k_{n}{+}m{+}1)$ such that
$(j_0,\cdots,j_m)\in\omega_m$ and $(k_0,\cdots,k_{n})\in\omega_n$.
}\end{Definition}\vspace{4mm}

\begin{Theorem}{\rm $(i_0,\cdots,i_n)\in\omega_n$ if and only if it is
the weight of a non-empty weight subgraph $G(i_0,\cdots,i_n)$.
}\end{Theorem}

Proof. We use induction on $n$. For $n=0,1$, the theorem is trivial.
Suppose the theorem holds for $n>1$. Then for $[a_{i,j}]\in
G(i_0,\cdots,i_n,i_{n+1})$, $i_{n+1}=\sum_{i=0}^n(1{-}a_{i,n+1})$,
\begin{eqnarray*}\left[\begin{array}{ccc}
a_{0,1}&\cdots&a_{0,n}\\
&\cdots&a_{1,n}\\
&&\cdots\\
&&a_{n-1,n}
\end{array}\right]\in G(i_0{-}a_{0,n+1},i_1{-}a_{1,n+1},\cdots,i_n{-}a_{n,n+1})
\end{eqnarray*}

By the induction hypothesis, for all
$0\!\leqslant\!k_0\!<\!\cdots\!<\!k_s\!\leqslant\!n$,
$(i_{k_0}{-}a_{k_0,n+1})\!+\!\cdots
\!+\!(i_{k_s}{-}a_{k_s,n+1})\geqslant 0\!+\!1\!+\!\cdots\! +\!s$.
So\\
\hspace*{30mm}$i_{k_0}\!+\!\cdots \!+\!i_{k_s}\geqslant
0\!+\!1\!+\!\cdots\! +\!s\!+\!\sum_{u=1}^sa_{k_u,n+1}\geqslant
0\!+\!1\!+\!\cdots\!
+\!s$,\\
\hspace*{5mm}$i_{k_0}\!+\!\cdots \!+\!i_{k_s}\!+\!i_{n+1}\geqslant
0\!+\!1\!+\!\cdots\!
+\!s\!+\!\sum_{u=1}^sa_{k_u,n+1}\!+\!\sum_{i=0}^n(1{-}a_{i,n+1})\geqslant
0\!+\!1\!+\!\cdots\! +\!s\!+\!(s{+}1),$\\
\hspace*{10mm}$i_0\!+\!\cdots\!i_n\!+\!i_{n+1}=i_0\!+\!\cdots\!i_n\!
+\!\sum_{i=0}^na_{i,n+1}\!+\!\sum_{i=0}^n(1{-}a_{i,n+1})
=0{+}1{+}\cdots{+}n{+}(n{+}1)$.\\
Thus, $(i_0,\cdots,i_n,i_{n+1})\in\omega_{n+1}$.

For $(i_0,\cdots,i_n,i_{n+1})\in\omega_{n+1}$, we will prove that
$G(i_0,\cdots,i_n,i_{n+1})$ is non-empty. By 1. of Theorem 3.1, we
may suppose
$i_0\!\leqslant\!i_1\!\leqslant\!\cdots\!\leqslant\!i_{n}\!\leqslant\!i_{n+1}$.
We use induction on $n{+}1{-}i_{n+1}$. If $i_{n+1}=n{+}1$, then
$G(i_0,\cdots,i_n,i_{n+1})=G(i_0,\cdots,i_n)$. The conclusion holds.
Suppose for $(j_0,\cdots,j_{n+1})\in\omega_{n+1}$ and
$j_{n+1}\!<\!n{+}1{-}s$, $G(j_0,\cdots,j_{n+1})$ is non-empty. Then
for $(i_0,\cdots,i_n,i_{n+1})\in\omega_{n+1}$ and
$i_{n+1}\!=\!n{+}1{-}s$, let $k$ be the biggest number such that
$i_0\!+\!i_1\!+\!\cdots\!+\!i_k=0\!+\!1\!+\!\cdots\!+\!k$. Then
$i_{k+1}\!>\!k{+}1$ and\\
\hspace*{10mm} $\alpha=(i_0,\cdots,i_k,i_{k+1}{-}1,\cdots,
i_{k+l}{-}1,i_{k+l+1},\cdots,i_{n},i_{n+1}{+}l)\in\omega_{n+1}$
($l=i_{k+1}{-}k{-}1$).\\
By the induction hypothesis, there exists $[a_{i,j}]\in G(\alpha)$.
We have $a_{k+1,n+1}=\cdots=a_{k+l,n+1}=0$. If not, then the product
$e=\prod e_{i,j}^{a_{i,j}}$ has a factor $e_{k+i,n+1}$. We may
suppose $i=1$. Then the weight of $e/e_{k+1,n+1}$ is
$(i_0,\cdots,i_k,i_{k+1}{-}2,\cdots,
i_{k+l}{-}1,i_{k+l+1},\cdots,i_{n},i_{n+1}{+}l{+}1)$. This
$(n{+}1)$-tuple is not in $\omega_{n+1}$. A contradiction! So
$ee_{k+1,n+1}\cdots e_{k+l,n+1}\in G(i_0,\cdots,i_n,i_{n+1})$. The
induction is complete. \hfill Q.E.D. \vspace{4mm}

\begin{Theorem}{\rm If $(k_0,{\cdots},k_{m+n+1})$  is reducible and is
the permutation of $(i_0,{\cdots},i_m,j_0{+}m{+}1,\\
\cdots,j_n{+}m{+}1)$ , then there is a GAD graph isomorphism\\
\hspace*{40mm}$G(k_0,\cdots,k_{m+n+1})= G(i_0,\cdots,i_{m}){\times}
G(j_0,\cdots,j_{n})$.\\ Specifically, for $m=0$ and $n=0$, we have\\
\hspace*{23mm}$G(j_0{+}1,\cdots,j_{k-1}{+}1,0,j_k{+}1,\cdots,
j_n{+}1)= G(j_0,\cdots,j_{k-1},j_k,\cdots, j_n),$\\
\hspace*{30mm}$G(i_0,\cdots,i_{k-1},m{+}1,i_k,\cdots,i_m)=
G(i_0,\cdots,i_{k-1},i_k,\cdots,i_m).$}\end{Theorem}

Proof. We consider the DGA $(R_n,\delta)$. For a reducible weight
$(k_0,\cdots,k_{m+n+1})$  that is a permutation of
$(i_0,\cdots,i_m,j_0{+}m{+}1,\cdots,j_n{+}m{+}1)$, there are
injective order preserving maps
$\sigma\colon\{0,\cdots,m\}\to\{0,\cdots,m{+}n{+}1\}$ and
$\tau\colon\{0,\cdots,n\}\to\{0,\cdots,m{+}n{+}1\}$ such that ${\rm
im}\sigma\cup{\rm im}\tau=\{0,{\cdots},m{+}n{+}1\}$, ${\rm
im}\sigma\cap{\rm im}\tau=\phi$. Then, the two maps induce two
algebra monomorphisms $\sigma_*\colon \Lambda(\frak A_{m+1})\to
\Lambda(\frak A_{m+n+1})$ and $\tau_*\colon \Lambda(\frak
A_{n+1})\to \Lambda(\frak A_{m+n+1})$ defined by
$\sigma_*(x_{s,t})=x_{\sigma(s),\sigma(t)}$ and
$\tau_*(x_{s,t})=x_{\tau(s),\tau(t)}$. $\sigma_*$ and $\tau_*$ are
not DGA homomorphisms. But let $c=\prod_{s,t} x_{\tau(s),\sigma(t)}$
(in any fixed order). Then, define linear map $\xi\colon
R(i_0,\cdots,i_{m})\otimes R(j_0,\cdots,j_{n})\to
R(k_0,\cdots,k_{m+n+1})$ by that $\xi(a\otimes
b)=\sigma_*(a)\tau_*(b)c$ for all $a\in R(i_0,\cdots,i_{m})$ and
$b\in R(j_0,\cdots,j_{n})$. It is obvious that $\delta c=0$ and
$\delta \Big(\sigma_*(a)\tau_*(b)c\Big)=\Big(\delta
\sigma_*(a)\Big)\tau_*(b) c+(-1)^{|a|}\sigma_*(a)\Big(\delta
\tau_*(b)\Big)c$. Thus, $\xi$ is a DGA monomorphism.

To prove that $\xi$ is an epimorphism, we need only show the two
free groups have the same dimension. By 1. of Theorem 3.1, we may
suppose
$(k_0,\cdots,k_{m+n+1})=(i_0,\cdots,i_m,j_0{+}m{+}1,\\\cdots,j_n{+}m{+}1)$.
For $[a_{i,j}]\in R(i_0,\cdots,i_m,j_0{+}m{+}1,\cdots,j_n{+}m{+}1)$,
\begin{eqnarray*}&&\quad\left[\begin{array}{ccc}
a_{0,1}&\cdots&a_{0,m+n}\\
&\cdots&a_{1,m+n}\\
&&\cdots\\
&&a_{m+n-1,m+n}
\end{array}\right]\\
&&\in
R(i_0{-}a_{0,m+n+1},\cdots,i_m{-}a_{m,m+n+1},j_0{-}a_{m+1,m+n+1}\cdots,j_{n-1}{-}a_{m+n-1,m+n+1}).
\end{eqnarray*}
Therefore, $i_0{-}a_{0,m{+}n{+}1}+\cdots
i_m{-}a_{m,m{+}n{+}1}\geqslant \frac{m(m{+}1)}{2}$ and so
$a_{k,m+n+1}=0$ for $k=0,\cdots,m$. This implies that
$R(i_0{-}a_{0,m+n+1},\cdots,i_m{-}a_{m,m+n+1},j_0{-}a_{m+1,m+n+1}\cdots,j_{n-1}{-}a_{m+n-1,m+n+1})$
is also reducible and
\begin{eqnarray*}&&\quad\left[\begin{array}{ccc}
a_{0,1}&\cdots&a_{0,m+n-1}\\
&\cdots&a_{1,m+n-1}\\
&&\cdots\\
&&a_{m+n-2,m+n-1}
\end{array}\right]\\
&&\in
R(i_0{-}a_{0,m+n},\cdots,i_m{-}a_{m,m+n},\cdots,j_{n-2}{-}a_{m+n-2,m+n}{-}a_{m+n-2,m+n+1}).
\end{eqnarray*}
For the same reason, $a_{k,m+n}=0$ for $k=0,\cdots,m$. Inductively,
we have that $a_{k,m+l+1}=0$ for $k=0,\cdots,m$ and $l=0,\cdots,n$.
That is,

\begin{eqnarray*}[a_{i,j}]=\left[\begin{array}{ccccccc}
a_{0,1}&\cdots&a_{0,m}&0&0&\cdots&0\\
&\cdots&\cdots&\cdots&\cdots&\cdots&\cdots\\
&&a_{m-1,m}&0&0&\cdots&0\\
&&&0&0&\cdots&0\\
&&&&a_{m+1,m+2}&\cdots&a_{m+1,m+n+1}\\
&&&&&\cdots&\cdots\\
&&&&&&a_{m+n,m+n+1}\end{array}\right].
\end{eqnarray*}
So dim$R(i_0,\cdots,i_m,j_0{+}m{+}1,\cdots,j_n{+}m{+}1)=$
dim$R(i_0,\cdots,i_m)\times$dim$R(j_0,\cdots,j_n)$. \hfill Q.E.D.
\vspace{4mm}

Notice that the cohomology classes with reducible weight may not be
decomposable. For example, the cohomology class of
$H^*(R(i_0,\cdots,i_n))$ with $(i_0,\cdots,i_n)$ a permutation of
$(0,1,\cdots,n)$ is not always decomposable, but all such weights
are reducible. \vspace{4mm}

\begin{Theorem}{\rm For every $\alpha\in\omega_n$,
$G(\alpha)$ is connected. Thus, $\omega_n$ is the set of connected
components of the graph $G_n=\Lambda(\frak A_{n+1})$. Moreover,  if
both $(i_0,\cdots,i_s,\cdots,i_t,\cdots,i_n)$ and
$(i_0,\cdots,i_s{+}1,\cdots,i_t{-}1,\cdots,i_n)$ are in $\omega_n$,
then there is at least one $e\in G(i_0,\cdots,i_n)$ such that
$ee_{s,t}\in G(i_0,\cdots,i_s{+}1,\cdots,i_t{-}1,\cdots,i_n)$.
}\end{Theorem}

Proof.  Firstly, we prove the second conclusion. By 1. of Theorem
3.1, we need only prove the case that if both $(i_0,i_1,\cdots,i_n)$
and $(i_0{+}1,i_1{-}1,\cdots,i_n)$ are in $\omega_n$, then there is
$e\in G(i_0,i_1,\cdots,i_n)$ such that $ee_{0,1}\in
G(i_0{+}1,i_1{-}1,\cdots,i_n)$. Suppose $i_0\!\leqslant\!i_1$. Let
$e=[a_{i,j}]\in G(i_0,i_1,\cdots,i_n)$. If $a_{0,1}=0$, then
$ee_{0,1}\in G(i_0{+}1,i_1{-}1,\cdots,i_n)$. If $a_{0,1}=1$, since
$\sum_{i=2}^n a_{0,i}<\sum_{i=2}^n a_{0,i}+1\leqslant\sum_{i=2}^n
a_{1,i}$, there is $2\!\leqslant\!i\!\leqslant\!n$ such that
$a_{1,i}=1$ and $a_{0,i}=0$. Let $e'=(e/e_{0,1}e_{1,i})e_{0,i}$,
then $e'\in G(i_0,i_1,\cdots,i_n)$ and $e'e_{0,1}\in
G(i_0{+}1,i_1{-}1,\cdots,i_n)$.  Suppose $i_0\!>\!i_1$. Firstly,
apply the conclusion to $(n{-}i_0,n{-}i_1,\cdots,n{-}i_n)$. Then
apply the duality isomorphism  $\vartheta$ in 4. of Theorem 3.1, we
prove the case for $i_0\!>\!i_1$.

Now we use induction on $n$ to prove that for every
$(i_0,\cdots,i_n)\in\omega_n$, $G(i_0,\cdots,i_n)$ is connected. For
$n=0,1,2$, it is a direct checking. Now suppose the conclusion holds
for $n\!>\!2$. Then for $(i_0,\cdots,i_n,i_{n+1})\in\omega_{n+1}$,
we have a disjoint union of sets
$$G(i_0,\cdots,i_n,i_{n+1})=\bigcup_{\varepsilon_0+\cdots+\varepsilon_n=n+1-i_{n+1}}
\;G(i_0{-}\varepsilon_0,\cdots,i_n{-}\varepsilon_{n})e^{\varepsilon_0}_{0,n+1}\cdots
e^{\varepsilon_n}_{n,n+1},$$ where the union is taken throughout all
$\varepsilon_s=0$ or $1$ such that
$(i_0{-}\varepsilon_0,\cdots,i_n{-}\varepsilon_{n})$ is a weight. By
the induction hypothesis, every
$G(i_0{-}\varepsilon_0,\cdots,i_n{-}\varepsilon_{n})$ is connected
in $G_n$ and so every
$G(i_0{-}\varepsilon_0,\cdots,i_n{-}\varepsilon_{n})e^{\varepsilon_0}_{0,n+1}\cdots
e^{\varepsilon_n}_{n,n+1}$ is connected in $G_{n+1}$. So we need
only prove that these different subgraphs are joined by paths in
$G_{n+1}$. Let $(i_0{-}\varepsilon_0,\cdots,i_n{-}\varepsilon_{n})$
and $(i_0{-}\varepsilon'_0,\cdots,i_n{-}\varepsilon'_{n})$ be
weights of two different components of the above disjoint union such
that $\varepsilon_s\!=\!\varepsilon'_t\!=\!0$,
$\varepsilon_t\!=\!\varepsilon'_s\!=\!1$ for some
$0\!\leqslant\!s\!<\!t\!\leqslant\!n$,
$\varepsilon_i\!=\!\varepsilon'_i$ for $i\neq s,t$. By the second
conclusion of the theorem, there is
$e\!\in\!G(i_0{-}\varepsilon_0,\cdots,i_n{-}\varepsilon_{n})$ such
that
$ee_{s,t}\!\in\!G(i_0{-}\varepsilon'_0,\cdots,i_n{-}\varepsilon'_{n})$.
Then,
\begin{eqnarray*}ee_{s,t}
e^{\varepsilon'_0}_{0,n+1}\cdots e^{\varepsilon'_n}_{n,n+1} &\in&
G(i_0{-}\varepsilon'_0,\cdots,i_n{-}\varepsilon'_{n})
e^{\varepsilon'_0}_{0,n+1}\cdots e^{\varepsilon'_n}_{n,n+1}\\
ee^{\varepsilon_0}_{0,n+1}\cdots e^{\varepsilon_n}_{n,n+1} &\in&
G(i_0{-}\varepsilon_0,\cdots,i_n{-}\varepsilon_{n})
e^{\varepsilon_0}_{0,n+1}\cdots e^{\varepsilon_n}_{n,n+1}
\end{eqnarray*}
But $[ee_{s,t} e^{\varepsilon'_0}_{0,n+1}\cdots
e^{\varepsilon'_n}_{n,n+1},ee^{\varepsilon_0}_{0,n+1}\cdots
e^{\varepsilon_n}_{n,n+1}]$ is an edge of $G_{n+1}$. Therefore, such
two components of the above disjoint union can be joined by paths in
$G_{n+1}$. This case implies that any two components of the above
disjoint union can be joined by paths in $G_{n+1}$. So
$G(i_0,\cdots,i_n,i_{n+1})$ is connected.\hfill Q.E.D.\vspace{4mm}

\begin{Theorem}{\rm The ranks of the weight subgraphs satisfy the
formula\\
\hspace*{30mm}$r(G(\cdots,j{+}1,\cdots,i{-}1,\cdots))-r(G(\cdots,j,\cdots,i,\cdots))=i{-}j{-}1$,\\
where the omitted parts of the two weights are the same. Thus, for a
weight $(i_0,i_1,\cdots,i_n)$ such that
$i_0\leqslant i_1\leqslant\cdots\leqslant i_n$,\\
\hspace*{40mm}$r(G(i_0,i_1,\cdots,i_n))=\sum_{k=0}^n
m_k+\sum_{i_k<k}(k-i_k),$\\
where $m_k=i_k{+}(i_k{+}1){+}\cdots{+}(k{-}1)$ if $i_k<k$ ;
$m_k=-i_k{-}(i_k{-}1){-}\cdots{-}(k{+}1)$ if $i_k>k$ ; $m_k=0$ if
$i_k=k$.}\end{Theorem}

Proof. Suppose $j$ is the $s$-th number and $i$ is the $t$-th number
of the weight $(\cdots j\cdots i\cdots)$. By the second conclusion
of Theorem 3.4, there is $a=[a_{u,v}]\!\in\!G(\cdots j\cdots
i\cdots)$ such that $a_{s,t}=0$ and $ae_{s,t}\!\in\!G(\cdots
j{+}1\cdots i{-}1\cdots)$. Let $b=[b_{u,v}]$ be the corresponding
triangular matrix for $aa_{s,t}$. The following lists all the
neighbors of $a$ and $b$.

(1) For $m<s$, $a_{m,s}=0$, $a_{m,t}=1$, let $a_m=[c_{u,v}]$ be
defined as follows. $c_{m,s}=1$, $c_{s,t}=1$, $c_{m,t}=0$,
$c_{u,v}=a_{u,v}$, otherwise. Then $a_m$ is a neighbor of $a$.

(2) For $s<n<t$, $a_{s,n}=1$, $a_{n,t}=1$, let $a_n=[c_{u,v}]$ be
defined as follows. $c_{s,n}=0$, $c_{n,t}=0$, $c_{s,t}=1$,
$c_{u,v}=a_{u,v}$, otherwise. Then $a_n$  is a neighbor of $a$.

(3) For $t<l$, $a_{s,l}=1$, $a_{t,l}=0$, let $a_l=[c_{u,v}]$ be
defined as follows. $c_{s,l}=0$, $c_{s,t}=1$, $c_{t,l}=1$,
$c_{u,v}=a_{u,v}$, otherwise. Then $a_l$  is a neighbor of $a$.

(4) For $m<s$, $b_{m,s}=1$, $b_{m,t}=0$, let $b_m=[c_{u,v}]$ be
defined as follows. $c_{m,s}=0$, $c_{s,t}=0$, $c_{m,t}=1$,
$c_{u,v}=b_{u,v}$, otherwise. Then $b_m$ is a neighbor of $b$.

(5) For $s<n<t$, $b_{s,n}=0$, $b_{n,t}=0$, let $b_n=[c_{u,v}]$ be
defined as follows. $c_{s,n}=1$, $c_{n,t}=1$, $c_{s,t}=0$,
$c_{u,v}=b_{u,v}$, otherwise. Then $b_n$ is a neighbor of $b$.

(6) For $t<l$, $b_{s,l}=0$, $b_{t,l}=1$, let $b_l=[c_{u,v}]$ be
defined as follows. $c_{s,l}=1$, $c_{s,t}=0$, $c_{t,l}=0$,
$c_{u,v}=b_{u,v}$, otherwise. Then $b_l$ is a neighbor of $b$.

Thus, $r(b)-r(a)=$ the number of $b_k$'s $-$ the number of $a_k$'s
$=i{-}j{-}1$.

Notice that $G(0,1,\cdots,n)$ has only one vertex, the triangular
matrix with all entries 0 (1 of $\Lambda(\frak A_{n+1})$). So
$r(G(0,1,\cdots,n))=0$. Then the second formula is an easy induction
on $r(G(i_0,\cdots,i_n))$.\hfill Q.E.D.

\end{document}